\documentclass[11pt,leqno]{article}

\usepackage{amsmath,amsfonts,amscd,amssymb,theorem}

\long\def\comment#1\endcomment{}

\makeatletter
\begingroup
\gdef\th@dotted{\normalfont\itshape
  \def\@begintheorem##1##2{%
        \item[\hskip\labelsep \theorem@headerfont ##1\ ##2.]}%
\def\@opargbegintheorem##1##2##3{%
   \item[\hskip\labelsep \theorem@headerfont ##1\ ##2\ (##3).]}}
\endgroup
\makeatother

\theoremstyle{dotted}

\newtheorem{theorem}{Theorem}[section]
\newtheorem{lemma}[theorem]{Lemma}

\newtheorem{prop}[theorem]{Proposition}
\newtheorem{corr}[theorem]{Corollary}

\makeatletter
\begingroup
\gdef\th@upshape{\normalfont
  \def\@begintheorem##1##2{%
        \item[\hskip\labelsep \theorem@headerfont ##1\ ##2.]}%
\def\@opargbegintheorem##1##2##3{%
   \item[\hskip\labelsep \theorem@headerfont ##1\ ##2\ (##3).]}}
\endgroup
\makeatother

\theoremstyle{upshape}

\newtheorem{defn}[theorem]{Definition}
\newtheorem{remark}[theorem]{Remark}

\makeatletter
\renewcommand{\subsection}{\@startsection{subsection}{2}{0pt}{-3ex
plus -1ex minus -0.2ex}{-2mm plus -0pt minus
-2pt}{\normalfont\bfseries}} \makeatother

\makeatletter
\@addtoreset{equation}{section}
\makeatother

\newcommand{\cntrct}                
{\hspace{2pt}\raisebox{1pt}{\text{$\lrcorner$}}\hspace{2pt}}

\newcommand{\proof}[1][Proof.]{\smallskip\noindent{\em #1}}
\def\endproof{\hfill\ensuremath{\square}\par\medskip}

\def\eqref#1{\thetag{\ref{#1}}}

\let\latexref=\ref
\def\ref#1{{\normalfont{\latexref{#1}}}}

\newcommand{\wt}{\widetilde}

\newcommand{\bekar}{\natural}
\newcommand{\wrth}{\textstyle\int}
\newcommand{\dg}{\dagger}

\newcommand{\lle}{<\!\!<}

\newcommand{\idot}{{\:\raisebox{1pt}{\text{\circle*{1.5}}}}}
\newcommand{\hdot}{{\:\raisebox{3pt}{\text{\circle*{1.5}}}}}
\newcommand{\Z}{{\mathbb Z}}

\newcommand{\eps}{\varepsilon}
\renewcommand{\phi}{\varphi}

\newcommand{\F}{{\sf F}}
\newcommand{\G}{{\sf G}}

\newcommand{\Ll}{{\sf L}}
\newcommand{\Rr}{{\sf R}}

\newcommand{\kk}[1]{{\sf K}_\idot\langle#1\rangle}
\newcommand{\kkk}[2]{{\sf K}_{#2}\langle#1\rangle}
\newcommand{\cc}[1]{{\sf C}_\idot\langle#1\rangle}
\newcommand{\ii}[1]{{\sf I}\langle#1\rangle}

\newcommand{\Fr}{{\sf Fr}}

\newcommand{\Hom}{\operatorname{Hom}}
\newcommand{\Ext}{\operatorname{Ext}}
\newcommand{\Tor}{\operatorname{Tor}}

\newcommand{\Ker}{\operatorname{Ker}}

\newcommand{\Fun}{\operatorname{Fun}}

\newcommand{\id}{\operatorname{\sf id}}
\newcommand{\Id}{\operatorname{\sf Id}}
\newcommand{\gr}{\operatorname{\sf gr}}

\newcommand{\Tot}{\operatorname{{\sf Tot}}}

\newcommand{\tw}{ {(1)} }

\newcommand{\D}{{\cal D}}
\newcommand{\DF}{{\cal D}{\sf F}}
\newcommand{\C}{{\cal C}}

\newcommand{\B}{{\sf B}}
\newcommand{\Bb}{{\cal B}}
\newcommand{\bB}{\overline{B}}

\newcommand{\Hh}{{\cal H}}

\newcommand{\m}{{\mathfrak m}}
\newcommand{\calo}{{\cal O}}

\newcommand{\Cycl}{\operatorname{Cycl}}
\newcommand{\Sets}{\operatorname{Sets}}
\newcommand{\Maps}{\operatorname{Maps}}

\newcommand{\mmod}{{\text{\rm -mod}}}
\newcommand{\bimod}{{\text{\rm -bimod}}}
\newcommand{\ppt}{{\sf pt}}

\newcommand{\Spec}{\operatorname{Spec}}
\newcommand{\Shv}{\operatorname{Shv}}

\newcommand{\cchar}{\operatorname{\sf char}}

\newcommand{\DT}{\operatorname{\sf DT}}

\title{Non-commutative Cartier operator and Hodge-to-de Rham
degeneration}

\author{D. Kaledin}

\begin{document}

\maketitle

\tableofcontents

\section*{Introduction}

When one thinks about differential forms and the de Rham complex in
the context of homological algebra and algebraic geometry, one
usually considers de Rham cohomology and treats it as one more
cohomology theory for algebraic varieties, similar to Betti
cohomology, \'etale cohomology, crystalline cohomology in positive
characteristic, and so on. However, there is an alternative point of
view which has been developping slowly over the years and recently
became quite prominent.

Since the pioneering paper \cite{HKR}, it has been known that a
differential form on the spectrum of a smooth commutative algebra
can be treated as its Hochschild homology class. In 1983, it was
discovered, simultaneously and independently in \cite{C},
\cite{FT1}, \cite{LQ}, that the de Rham differential also has such
an algebraic interpretation. Moreover, and most surprisingly, the
whole formalism makes sense for an associative, but not necessarily
commutative algebra $A$. The de Rham cohomology appears in this
setting as {\em additive $K$-theory} (\cite{FT}) or {\em cyclic
homology} (\cite{C}). Hochschild and cyclic homology are related by
a spectral sequence generalizing the commutative Hodge to de Rham
spectral sequence. For an introduction to the subject, we refer the
reader to \cite{L} (and also to \cite{FT}, which contains much
material not covered in \cite{L}).

A natural next step in this direction would be to study the Hodge
theory -- since the standard spectral sequence relating Hodge
cohomology and de Rham cohomology makes sense for a non-commutative
algebra $A$, under what assumptions does this spectral sequence
degenerate? Recently there has been much interest in this topic --
most notably, in the work of M. Kontsevich, who has stated and
popularized a conjectural non-commutative analog on the standard
degeneration theorem for smooth projective algebraic varieties
(\cite{K1}, \cite{K2}). He also gave some very beautiful
applications, and formulated a set of very precise finiteness
assumptions which are necessary for the degeneration and play the
role of both properness and smoothness in the non-commutative
setting.

To the best of our knowledge, so far there has been very little
progress in proving Kontsevich's conjecture. This is perhaps not
very surprising, since in the non-commutative setting, the usual
analytic approach to Hodge theory makes no sense. Fortunately, there
is an alternative algebraic approach. It was discovered by
P. Deligne and L. Illusie \cite{DL} back in 1987, a few years after
the cyclic homology was discovered. The argument uses reduction to
positive characteristic and the Cartier isomorphism between de Rham
cohomology spaces and spaces of forms for smooth varieties over a
field of positive characteristic.

\medskip

The goal of the present paper is to apply the Deligne-Illusie method
in the non-commuative setting and to prove that the non-commutative
Hodge to de Rham spectral sequence also degenerates, provided some
natural finiteness conditions are met. We do this by giving a
version of the Cartier isomorphism valid without any assumptions of
commutativity. The final result is a direct generalization of
\cite{DL}. In particular, it reduces to \cite{DL} for smooth proper
commutative algebraic varieties (except that our general
construction of the Cartier isomorphism is somewhat more canonical
and does not use any explicit formulas at all).

\medskip

Our approach generally follows that of Kontsevich, but differs from
his in one important point. Kontsevich gave his conjecture in
characteristic $0$, and he used the language of
$A_\infty$-algebras. Our method is reduction to positive
characteristic, and the notion of $A_\infty$-algebra is not very
convenient there (because it uses the ``naive'', not simplicial
tensor structure on the category of complexes of vector spaces). The
same goes for differential-graded algebras. Instead, we give and
prove the statement for sheaves of usual associative algebras over
an arbitrary site. This includes both the case of usual algebraic
varieties, possibly considered with some non-commutative enhancement
of the structure sheaf, and the case of simplicial or cosimplicial
algebras. Most likely, the proper generality for the theory is in
any case that of a triangulated category equipped with some enriched
structure; however, this should be the topic of further research.
In one respect at least, $A_\infty$ approach is definitely better:
when working with sheaves, in order to be able to reduce the problem
in characteristic $0$ to a problem in positive characteristic, we
have to impose additional assumptions -- for instance, to assume
that the algebra sheaf in question is Noetherian. This is probably
too strong, and indeed, the finiteness assumptions needed anyway for
degeneration should also be sufficient for the reduction problem. In
a subsequent paper, we will investigate the relation between our
approach and the $A_\infty$ statement, and we hope to be able to
prove Kontsevich's conjecture in full generality.

We note that even in the commutative case, one might get a statement
by our method for at least some proper algebraic varieties that are
not smooth (the Hochschild homology in this case should be
understood as the Hochschild homology of the category of perfect
complexes of coherent sheaves). We do not know whether this has any
geometric significance. Also, it is known that the Deligne-Illusie
method allows one to prove several strong vanishing theorems of the
Kodaira type, see \cite{EV}; we do not know whether our
non-commutative version can be used in a similar way.

Another intriguing observation concernes the non-commutative Cartier
isomorphism. In a sense, and this may be made quite precise, the
Cartier isomorphism is only the visible part of an iceberg, which is
the action of the Frobenius map on crystalline cohomology. In the
non-commutative setting, there is certainly no Frobenius
map. However, our results suggest that the Frobenius action on
cohomology still exists. Among the many and varied applications of
this Frobenius action, some might also make sense in the
non-commutative world.

\subsection*{Acknowledgements} I am very grateful to M. Kontsevich
for his enthusiasm for this project and many valuable
discussions. His talk \cite{K1} stimulated my interest in the
problem. Another source of the present project is a question posed
by N. Markarian. Discussions with L. Katzarkov were also very
helpful. I am also grateful to A. Beilinson, R. Bezrukavnikov,
F. Bogomolov, A. Bondal, A. Braverman, B. Feigin, V. Ginzburg,
D. Kazhdan, A. Kuznetsov, S. Loktev, N. Markarian and D. Orlov. Part
of this work was done during a visit to the Hebrew University of
Jerusalem.

\section{A general overview.}

\subsection{Recollection on the commutative case.}

We start by briefly recalling the setup and the method of
\cite{DL}. Let $X$ be a smooth algebraic variety over a field
$k$. Then the cotangent sheaf $\Omega(X/k)$ is flat, so that we have
the sheaves $\Omega^\hdot$ of higher-degree differential forms, and
the de Rham differential $d:\Omega^\hdot(X) \to \Omega^{\hdot +
1}(X)$. One considers the de Rham complex $\langle \Omega^\hdot, d
\rangle$, and one equips it with the so-called {\em stupid
filtration}:
$$
F^i\Omega^j(X) = \begin{cases}
\Omega^j(X), &\qquad j \geq i,\\
0, &\qquad \text{otherwise}.
\end{cases}
$$
It is known that if $X$ is projective over $k$, and $\cchar k = 0$,
then the spectral sequence associated to this stupid filtration
degenerates, so that the de Rham cohomolgy groups $H^i_{DR}(X)$ of
the variety $X$ are isomorphic to the Hodge cohomology groups
$\bigoplus_{p+q=i}H^{p,q}(X)$, $H^{p,q}(X)=(X,\Omega^p(X))$.

The condition $\cchar k = 0$ is essential: a long time ago
D. Mumford constructed a counterexample, already in $\dim 2$, which
shows that the Hodge-to-de Rham spectral sequence does {\em not}
always degenerate in positive characteristic.

Nevertheless, Deligne and Illusie give a proof of the Hodge-to-de
Rham degeneration by using reduction to positive characteristic; at
the same time, they gain an understanding of the situation in
$\cchar p$ and show that the degeneration does hold when one
additional natural condition is imposed.

Namely, assume that $k$ is a perfect field, $\cchar k = p$ is a
positive prime, and assume that $\dim X < p$. Then the de Rham
complex $\Omega^\hdot(X)$ is perfectly well-defined, but it is not
exact even analytically, in a formal neighborhood of a point: it has
been shown by P. Cartier that there exists a canonical isomorphism
$$
C:\Hh^\hdot_{DR}(X) \cong \Omega^\hdot(X^\tw)
$$
between the de Rham cohomology sheaves $\Hh^\hdot_{DR}(X)$ of the
variety $X$ and the sheaves $\Omega^\hdot(X^\tw)$ of forms on the
Frobenius twist $X^\tw$ (under our assumptions, $X^\tw$ is
isomorphic to $X$ equipped with the structure sheaf $\calo_X^p$, the
subsheaf of $p$-th powers in $\calo_X$). In particular, the kernel
of the de Rham differential $d:\calo_X \to \Omega^1(X)$ is precisely
the subsheaf $\calo_X^p$, so that the de Rham differential is
$\calo_X^p$-linear (and Zariski topology is fine enough for all
computations with the de Rham complex). The Cartier isomorphism
$\Hh^0_{DR}(X) \cong \Omega_{X^\tw}$ is inverse to the Frobenius
map.

The crucial part of \cite{DL} is concerned with the following
question: when does the Cartier isomorphism $C$ -- or rather, the
inverse map $C^{-1}$ -- extend to a map of complexes? This can be
rephrased as follows. Consider the {\em canonical} filtration $\tau$
on the de Rham complex $\Omega^\hdot(X)$ -- namely, let
$$
\tau_{\leq i}\Omega^j(X) = \begin{cases}
\Omega^j(X), &\qquad j < i,\\
\Ker d, &\qquad j=i,\\
0, &\qquad \text{otherwise}.
\end{cases}
$$
Then the associated graded quotient with respect to this canonical
filtration is naturally quasiisomorphic to the sum
$\bigoplus_i\Hh_{DR}(X)[-i] \cong
\bigoplus_i\Omega^i(X^\tw)[-i]$. We want to know when the de Rham
complex is quasiisomorphic to its associated graded quotient -- in
other words, when the canonical filtration splits.

The main result of \cite{DL} claims that -- and this is a purely
local statement -- this happens if and only if the variety $X$
admits a lifting to a variety smooth over the second Witt vectors
ring $W_2(k)$ (and moreover, these splittings satisfying some
natural conditions are in one-to-one correspondence with liftings of
$X$ to $W_2(k)$).

The rest of the proof is a surprisingly easy corollary of this basic
fact. The canonical filtration induces a filtration on de Rham
cohomology groups $H^\hdot_{DR}(X)$ known as {\em conjugate
filtration}, and a spectral sequence called conjugate spectral
sequence -- its first term consists of the same Hodge cohomology
groups $H^{p,q}(X) \cong H^{p,q}(X^\tw)$ as for the Hodge spectral
sequence, but the conjugate and the Hodge filtrations go in the
opposite direction. By the cited main result of \cite{DL}, if $X$
can be lifted to $W_2(k)$, then the conjugate spectral sequence
degenerates. But if $X$ is projective -- or in fact just proper --
over $k$, all these spectral sequences consist of {\em
finite-dimensional} $k$-vector spaces. Thus we have two spectral
sequences of finite-dimensional vector spaces with the same term
$E^1$ and the same term $E^\infty$: for dimension reasons, if one
degenerates, the other must degenerate, too.

\subsection{Cyclic approach.}

Assume now that $X = \Spec A$ is affine. The {\em Hochschild
homology} $HH_\idot(A)$ of the $k$-algebra $A$ is by definition the
$\Tor$-groups of the diagonal $A$-bimodule with itself: we have
$$
HH_\idot(A) = \Tor^\hdot_{A^{opp} \otimes A}(A,A).
$$
It has been established in \cite{HKR} that under the assumptions
above -- $X$ smooth and $\cchar k = 0$ -- we have $HH_l(A) \cong
\Omega^l(A)$ for any $l \geq 0$. In positive characteristic, this is
also true provided $l < \cchar k$ -- or always, if one requires
$\dim X < \cchar k$.

The Hochschild homology of any algebra can be computed by the bar
resolution -- this results in the well-known standard complex
consisting of tensor powers $A^{\otimes l}$, with a certain explicit
differential $b:A^{\otimes l+1} \to A^{\otimes l}$. The cyclic
homology $HC_\idot(A)$ has several equivalent definitions; by one of
them, $HC_\idot(A)$ is the total homology of a periodic bicomplex
\begin{equation}\label{connes}
\begin{CD}
@. @. @. @. A\\
@. @. @. @. @AA{b}A\\
@. @. @. A @>{B}>> A^{\otimes 2}\\
@. @. @. @AA{b}A @AA{b}A\\
@. @. A @>{B}>> A^{\otimes 2} @>{B}>> A^{\otimes 3}\\
@. @. @AA{b}A @AA{b}A @AA{b}A\\
@. A @>{B}>> A^{\otimes 2} @>{B}>> 
A^{\otimes 3} @>{B}>> A^{\otimes 4}
\end{CD}
\end{equation}
Here $b$ is the Hochschild differential, and $B$ is the new
differential introduced by A. Connes. In our situation, if one takes
the vertical cohomology, one obtains the complex
\begin{equation}\label{conn.dr}
\begin{CD}
@. @. @. @. A\\
@. @. @. @. @AA{0}A\\
@. @. @. A @>{d}>> \Omega^1(A)\\
@. @. @. @AA{0}A @AA{0}A\\
@. @. A @>{d}>> \Omega^1(A) @>{d}>> \Omega^2(A)\\
@. @. @AA{0}A @AA{0}A @AA{0}A\\
@. A @>{d}>> \Omega^1(A) @>{d}>> \Omega^2(A) @>{d}>> \Omega^3(A)
\end{CD}
\end{equation}
The Connes differential $B$ becomes the usual de Rham differential,
and the rows of this bicomplex are all truncation of the de Rham
complex.

Both $HH_\idot(A)$ and $HC_\idot(A)$ are defined for an arbitrary
associative algebra $A$, and so is the bicomplex \eqref{connes} --
but this is as far as one gets in the general situation: to split
\eqref{conn.dr} nicely into separate rows and consider them
separately, one needs to know that the algebra $A$ is
commutative. Nevertheless, for some applications this is not
needed. In particular, the stupid filtration on the de Rham complex
can be understood as the stupid filtration on the bicomplex
\eqref{connes} taken ``in horizontal direction'', and this makes
perfect sense in full generality (to avoid confusion, we note that
the modules of differential forms $\Omega^i(A)$ appear in the
Hochschild homology theory with ``wrong'' degrees -- $\Omega^l(A)$
appears in degree $-l$, not $l$ -- and as the result of this, stupid
a.k.a. Hodge filtration becomes increasing, not decreasing). Thus
one gets an analog of the Hodge-to-de Rham spectral sequence and can
pose the question of its degeneration (we note that this
``Hochschild-to-cyclic'' spectral sequence has been very prominent
in cyclic homology studies from the very beginning).

Moreover, in characteristic $\cchar k = p > 0$, one can define the
{\em conjugate} filtration on \eqref{conn.dr} as the canonical
filtration ``in horizontal direction'', and repeat the argument of
\cite{DL} (with some minor modifications). The end result is the
same -- one shows that the Hodge-to-de Rham spectral sequence
degenerates, once some natural conditions are met.

\subsection{Overview of the present paper.}

In a nutshell, what we do in the present paper is this: we note that
the conjugate filtration on \eqref{conn.dr}, known in the
commutative case, can be defined in an invariant way which makes
sense for a arbitrary associative algebra or sheaf of algebras over
the base. Moreover, we identify the associated graded quotient of
this filtration with the Hochschild homology -- thus obtaining a
version of the Cartier isomorphism valid in the non-commutative
setting (or rather, as in \cite{DL}, we in fact construct the
inverse map $C^{-1}$).

As an example, one can consider what in happens in degree $0$. Here
we have $HH_0(A) = HC_0(A) = A/[A,A]$, the quotient of $A$ by the
subspace spanned by commutator expressions $[a,b]=ab-ba$, $a,b, \in
A$. It is a pleasant exercize to check that, even for a
non-commutative $A$, the ``Frobenius'' map $x \mapsto x^p$ descends
to a well-defined and additive map $C^{-1}:A/[A,A] \to A/[A,A]$

It turns out that to generalize it to all degrees, it is not convenient
to use the complex \eqref{connes}, and it is better to use the
second standard complex for $HC_\idot(A)$, which is the periodic
complex
\begin{equation}\label{conn.2}
\begin{CD}
@>>> A @>>> A @>>> A @>>> A\\
@. @AA{b'}A @AA{b}A @AA{b'}A @AA{b}A\\
@>>> A^{\otimes 2} @>>> A^{\otimes 2} @>>> A^{\otimes 2} @>>>
A^{\otimes 2}\\
@. @AA{b'}A @AA{b}A @AA{b'}A @AA{b}A\\
@>>> A^{\otimes 3} @>>> A^{\otimes 3} @>>> A^{\otimes 3} @>>>
A^{\otimes 3}\\
@. @AA{b'}A @AA{b}A @AA{b'}A @AA{b}A\\
@>>> A^{\otimes 4} @>>> A^{\otimes 4} @>>> A^{\otimes 4} @>>>
A^{\otimes 4}
\end{CD}
\end{equation}
Here $b$ is the Hochschild differential, $b'$ is the contractible
differential in the bar resolution, and the horizotal differential
in $l$-th row is the same as in the standard periodic complex which
computes the homology $H_\idot(\Z/l\Z,A^{\otimes l})$ of the cyclic
group $\Z/l\Z$ with coefficients in $A^{\otimes l}$ equipped with
the natural action. Moreover, we note that the standard bar complex
which computes Hochschild homology can be modified in such a way
that instead of tensor powers $A^{\otimes l}$, it only contains
tensor powers $A^{\otimes pl}$ for some fixed integer $p$. This is
easy to see using the $\Tor$-interpretation of Hochschild homology:
since $A = A \otimes_{A^{opp} \otimes A} A$, we have
$$
HH_\idot(A) = \Tor^\hdot_{A^{opp} \otimes A}(A \otimes_{A^{opp}
  \otimes A}  A \otimes_{A^{opp} \otimes A} \cdots \otimes_{A^{opp}
  \otimes A} A,A),
$$
with an arbitrary fixed number of multiples on the left-hand
side. If one uses the bar resolution and multiplies it with itself
$p$ times in the correct way (see Lemma~\ref{bounded.dg}), the
result is a complex computing $HH_\idot(A)$ and consisting of
$A^{\otimes pl}$, $l \geq 0$.

Once one develops a version of \eqref{conn.2} for this modified
Hochschild complex, the Cartier isomorphism in degree $\geq 1$ boils
down essentially to the following linear-algebraic fact (this is
Lemma~\ref{otimesp}).
\begin{itemize}
\item For any vector space $V$ over a field $k$ of characteristic $p
  > 0$, the homology group $H_l(\Z/p\Z,V^{\otimes p})$ is
  canonically isomorphic to $V$ for every $l \geq 1$.
\end{itemize}
Unfortunately, we were not able to find an explicit complex to play
the role of \eqref{conn.2}; to overcome this difficulty, we use the
more invariant technique of homology of small categories. The
relevant categories are the Connes cyclic category $\Lambda$ and its
$p$-fold cover $\Lambda_p$, slightly less-known but also very
standard (see e.g. \cite{GB}). The invariant definition of the
conjugate filtration is obtained through considering the natural
functor $\Lambda_p \to \Lambda$ -- this functor is a fibration,
whose fiber is the groupoid with one object and automorphism group
$\Z/p\Z$, and the definition of the conjugate filtration that we use
is obtained by considering the homology of this group (it turns out
that the right thing to do is to consider the {\em even} terms of
the canonical filtration on the standard periodic complex).

We then study the question of the splitting of the conjugate
filtration. The end result, somewhat surprisingly, is again just the
same as in \cite{DL} -- the filtration splits if and only if $A$ can
be lifted to a flat algebra over the second Witt vector ring
$W_2(k)$. Since there is no natural notion of dimension for general
associative algebras, one might expect that the condition $\dim X <
\cchar k$ would disappear from the picture completely. This is not
what happens, though: our method only gives splitting only ``in
degrees between $0$ and $2p$'', and there reasons to believe that
this is not a merely technical limitation (see, in particular,
Remark~\ref{steen}). Nevertheless, since to deduce Hodge-to-de Rham
degeneration we are allowed to consider only generic primes, this
does not interfere with the proof (just as in \cite{DL}).

\bigskip

Let us now describe briefly the organization of the paper. A rather
long Section~\ref{intro} contains all the preliminary facts that we
need on homology of small categories in general and cyclic categories
$\Lambda_p$ in particular; nothing here is new, and everything is
contained in some form in \cite{L}. However, we do need to recall
these things in some detail to set up the notation etc. As the
result, the present paper is self-contained to a large degree. The
last part of Section~\ref{intro}, Subsection~\ref{para.sub}, is
taken with Proposition~\ref{i.p} which essentially explains how to
compute the cyclic homology of an algebra $A$ by means of a complex
only containing tensor powers $A^{\otimes pl}$ (although it is
formulated in the invariant language of small categories and cyclic
objects, and the statement sounds quite differently).

In Section~\ref{filtr}, we introduce the Hodge and the conjugate
filtration on cyclic homology of the arbitrary associative algebra
$A$, and we indetify the associated graded quotients of the
conjugate filtration with the Hochschild homology. Then in
Section~\ref{horror}, we study the extension data between the
associated graded pieces of the conjugate filtration. It turns out
that this is quite complicated technically, and moreover, that we
need to exploit some additional symmetry that we call
polycyclic. Partially, this technical complexity reflects the fact
that the resulting description only holds in degrees between $0$ and
$2p$ (while in higher degrees one encounters Steenrod cohomological
operations, and understanding the picture is probably far outside of
what our approach allows).

We note that, for better or for worse, we have tried to separate
consistently ``linear-algebraic'' facts, which make sense for any
cyclic object, and things that involve tensor products. Therefore in
Sections~\ref{filtr} and \ref{horror} we work in more generality
than strictly needed for the computation of $HC_\idot(A)$ for an
associative algebra $A$. Associative algebras {\em per se} and their
homology only appear in Section~\ref{alg}; with all the work done in
Section~\ref{filtr} and Section~\ref{horror}, in
Subsection~\ref{alg.def} it remains just to define a cyclic object
$A_{\#}$ associated to an algebra $A$ and to check that it satisfies
the assumptions of previous Sections (this is very easy). Then in
Subsection~\ref{spl.sub}, we study the splitting of the conjugate
filtration for an algebra $A$. With all the technical unpleasantness
already done in Section~\ref{horror}, this boils down to a nice and
compact criterion Lemma~\ref{dag.dege}; in Lemma~\ref{formm}, we
show by a reasonably simple and conceptual argument that the
condition of Lemma~\ref{dag.dege} is satisfied if $A$ admits a
lifting to $W_2(k)$.

Finally, in Section~\ref{main} we collect all of the above to prove
Hodge-to-de Rham degeneration. The statement in positive
characteristic, Theorem~\ref{char.p.main}, is a precise analog of
the corresponding commutative statement in \cite{DL}. Unfortunately,
the statement in $\cchar 0$, while also exactly the same as in the
commutative case, is not satisfactory: we need to require our
algebra to be Noetherian, which is too strong a condition in the
non-commutative case. It should be possible to drop this requirement
by using $A_\infty$ methods; we plan to return to this in the
future.

\section{Cyclic objects.}\label{intro}

In this section, we recall standard basic facts about cyclic objects
and cyclic categories needed in the rest of the paper; we mostly
follow \cite{FT}, see also \cite[Chapter 6]{L}. 

\setcounter{subsection}{-1}
\subsection{Notation.}
Throughout the paper, we will need to work with filtered complexes
in many places. Our conventions will be like this: if $\F^\hdot$ is
a filtered complex equipped with an increasing filtration $W_\idot$,
then $W_{[n,m]}\F^\hdot$ is the cone of the natural map
$W_{n-1}\F^\hdot \to W_m\F^\hdot$, and $\gr^W_n\F^\hdot$ is
$W_{[n,n]}\F^\hdot$. For decreasing filtrations, we use dual
conventions. The canonical filtration on a complex in a derived
category is denoted by $\tau_{\leq n}$; $\tau_{[n,m]}$ is the cone
of the map $\tau_{\leq (n-1)} \to \tau_m$. Sometimes we use
homological degrees $\F_\idot$ instead of cohomological degrees
$\F^\hdot$ -- we treat them as interchangeable notation for the same
thing, related by $\F^l = \F_{-l}$, $l \in \Z$.

\subsection{Homological preliminaries.}\label{prelim}
We begin with homological generalities. Let $\C$ be a Grothendieck
abelian category -- or, in Grothendieck's language \cite{G}, an
abelian category which additionally satisfies $AB3$, $AB4$, $AB3'$,
$AB4'$ (in applications, $\C$ will be the category of $B$-modules
for a commutative ring $B$, or, more generally, the category of
sheaves of $\Bb$-modules for a sheaf of commutative rings $\Bb$ on
some site).  For any small category $\F$, denote by $\Fun(\F,\C)$
the abelian category of covariant functors from $\F$ to $\C$, and
denote by $\D^-(\F,\C)$ the derived category of complexes in
$\Fun(\F,\C)$ bounded from above. The category $\Fun(\F,\C)$ is also
a Grothendieck abelian category.

For any two small categories $\F$, $\G$ and a functor $\sigma:\F \to
\G$, denote by $\sigma^*:\Fun(\G,\C) \to \Fun(\F,\C)$ the natural
restriction functor given by
$$
\sigma^*(E) = E \circ \sigma \in \Fun(\F,\C), \qquad E \in \Fun(\G,\C).
$$
Since $\C$ is a Grothendieck category, $\sigma^*$ has right and left
adjoint functors called {\em Kahn extensions}. Denote the left and
right Kahn extensions by $\sigma_!:\Fun(\F,\C) \to \Fun(\G,\C)$ and
$\sigma_*:\Fun(\F,\C) \to \Fun(\G,\C)$, and denote their derived
functors by $\Ll^\hdot\sigma_!:\D^-(\F,\C) \to \D^-(\G,\C)$ and
$\Rr^\hdot\sigma_!:\D^-(\F,\C) \to \D^-(\G,\C)$.

In the particular case when $\F=\ppt$ is the category with one
object and one morphism, we have $\Fun(\F,\C) \cong \C$, and
specifying a functor $\sigma:\F \to \G$ is the same as specifying an
object $[f] \in \G$. For any object $A \in C$, we denote
$$
\begin{aligned}
A_{[f]} &= \sigma_!A,\\
A_{[f]}^o &= \sigma_*A.
\end{aligned}
$$
If $\C$ is the category of abelian groups, and $A=\Z$ is the group
$\Z$, for any $g \in \G$ the group $\Z_{[f]}([g])$ is naturally
identified with the free abelian group spanned by the set
$\Hom_{\G}([f],[g])$, and the group $\Z_{[f]}^o([g])$ is naturally
identified with the abelian group of $\Z$-valued functions on the
set $\Hom_{\G}([g],[f])$.

In the particular case when $\G=\ppt$, there exists only one
tautological functor $\sigma:\F \to \G$. For any $A \in \C \cong
\Fun(\G,\C)$, we denote
$$
A_{\F} = \sigma^*A,
$$
or simply $A \in \Fun(\F,\C)$ when there is no danger of confusion.
For any $E \in \Fun(\F,\C)$ we denote $\Ll^\hdot\sigma_!(E)$ by
$H_\idot(\F,E)$ and we call it the {\em homology of the category
$\F$ with coefficients in $E \in \Fun(\F,\C)$}. Analogously, we
denote $\Rr^\hdot\sigma_*E = H^\hdot(\F,E)$, and we call the {\em
cohomology of the category $\F$ with coefficients in $E$}.  We note
that for any $\rho:\F \to \G$, and any $E \in \Fun(\G,\C)$, we have
by adjunction a canonical map
$$
\sigma_!:H_\idot(\F,\sigma^*\C) \to H_\idot(\G,C).
$$
If $\C = B\mmod$ is the category of modules over a commutative ring
$B$, we have by adjunction
$$
H^\hdot(\F,B_\F) = \Ext^\hdot(B_\F,B_\F),
$$
so that the cohomology $H^\hdot(\F,B_\F) = H^\hdot(\F,B)$ with
coefficients in the constant functor $B \in \Fun(\F,B\mmod)$ carries
a natural algebra structure. For any $B$-linear Grothendieck
category $\C$ and any object $E \in \Fun(\F,\C)$, the homology
$H_\idot(\F,E)$ is naturally a module over the algebra
$H^\hdot(\F,B_\F)$. For any $\F$, $\G$ and a functor $\rho:\F \to
\G$, the restriction functor $\rho^*$ induces an algebra map
$\rho^*:H^\hdot(\G,B) \to H^\hdot(\F,B)$. The maps $\rho^*$ and
$\rho_!$ are, as usual, related by the projection formula: we have
$$
\sigma_!(\sigma^*(\alpha) \cdot e) = \alpha \cdot \sigma_!(e)
$$
for any $\alpha \in H^\hdot(\G,B)$, $e \in H_\idot(\F,\rho^*E)$.

We will also need to use the notion of a fibered category
(\cite{SGA}); let us briefly recall it. For any functor $\sigma:\F \to
\G$ between small categories $\F$ and $\G$ and any object $[a] \in
\G$, denote by $\F_{[a]}$ the subcategory in $\F$ of all objects
$[a'] \in \F$ such that $\sigma([a'])=[a]$ and all morphisms $f$ such
that $\sigma(f)=\id_{[a]}$. This is called the {\em fiber of the
functor $\sigma:\F \to \G$ over $[a] \in \G$}. For any $[a],[b] \in
\G$, any $f:[a] \to [b]$ and any $[a'] \in \F_{[a]}$, $[b'] \in
\F_{[b]}$, denote by $\F_f([a'],[b'])$ the set of all maps $f':[a']
\to [b']$ such that $\sigma(f')=f$.

\begin{defn}
The functor $\sigma:\F \to \G$ is called a {\em fibration}, and $\F$
is called a {\em fibered category over $\G$}, if for any $[a],[b]
\in \G$, $f:[a] \to [b]$ and $[b'] \in \F_{[b]}$ the functor
$$
[a'] \mapsto \F_f([a'] \to [b'])
$$
from $\F_{[a]}$ to the category of sets is representable by an
object $f^*([b]) \in \F_{[a]}$.
\end{defn}

Passing to the opposite categories gives the dual notion of a
cofibration and a cofibered category. A functor which is both a
fibration and a cofibration is called a {\em bifibration}.

An example of a cofibered category is a so-called {\em discrete
cofibration}: for any functor $F$ from $\G$ to sets, pairs $\langle
[a] \in \G, a \in F([a]) \rangle$ form a small category $\F$
naturally fibered over $\G$ (the fiber $\F_{[a]}$ is the set
$F([a])$ cosidered as a category with no non-identical
morphisms). The category $\F$ is called the {\em total space} of the
functor $F$ and denoted $\Tot(\G,F)$. It is easy to check that every
cofibration whose fibers have no non-identical morphisms is obtained
in this way. An opposite case is a cofibration $\F \to \G$ whose
fibers are categories with one object. For example, every functor
$F$ from $\G$ to the category of groups defines a category $\F$
cofibered over $\G$ whose fiber $\F_{[a]}$ is a groupoid with one
object with automorphism group $F([a])$. By abuse of notation, in
this case we will also call $\F$ the total space of the functor $F$
and denote it by $\Tot(\G,F)$. More generally, a cofibered category
over $\G$ whose fibers are groupoids is usually called a {\em gerbe}
over $\G$. Gerbes of the form $\Tot(\G,F)$ are said to be {\em
split}; there are gerbes which are not of this type (we will see one
example below in Remark~\ref{nonspl}).

We will need the following standard properties of cofibrations (and
the dual properties of fibrations, which are analogous and left to
the reader).

\begin{prop}\label{fibr}
Let $\sigma:\F \to \G$ is a cofibration.
\begin{enumerate}
\item For any small category $\G'$ and a functor $\rho:\G' \to \G$,
  the pullback $\sigma':\G' \times_{\G} \F \to \G'$ is also a
  cofibration, and we have a natural base change isomorphism
\begin{equation}\label{base.change}
\Ll^\hdot\sigma_! \circ \rho^* \cong \rho^* \circ \Ll^\hdot\sigma_!.
\end{equation}
\item For any Grothendieck category $\C$ and any object $E \in
  \Fun(\G,\C)$, we have a natural isomorphism
\begin{equation}\label{proj.formula}
\Ll^\hdot\sigma_!\sigma^*E \cong E \otimes \Ll^\hdot\sigma_!\sigma^*\Z_{\G},
\end{equation}
and this isomorphism is functorial in $E$.
\end{enumerate}
\end{prop}

\proof[Sketch of a proof.] To check that $\sigma':\G' \times_{\G} \F
\to \G'$ is indeed a cofibration is an elementary exercize left to
the reader. Moreover, for any object $[a] \in \G'$ we have natural
isomorphism $\pi:\sigma'_{[a']} \cong \sigma_{\pi([a'])}$. Now, by the
definition of a cofibration, for any objects $M \in \C$, $[a] \in
\G$, we have
$$
\sigma^*M^o_{[a]} \cong \lim_{\gets}M^o_{[b]} \cong \iota_*(M_{\sigma_{[a]}}),
$$
where the limit is taken over $[b] \in \sigma_{[a]} \subset \F$, and
$\iota:\sigma_{[a]} \to \F$ is the embedding of the fiber
$\sigma_{[a]}$. Therefore for any $E \in \Fun(\F,\C)$ and any $[a] \in
\G$ we have
$$
\Ll^\hdot\sigma_!(E)([a]) \cong H_\idot(\sigma_{[a]},E),
$$
and this is functorial in $E$ and in $[a]$. This immediately gives
both the base change isomorphism \eqref{base.change} and the
projection formula \eqref{proj.formula}.
\endproof

Finally, in Section~\ref{filtr} we will need to use a bar resolution
several times; we recall the setup. Let $F:\C \to \C$ be any exact
functor from an abelian category $\C$ to itself. Assume that $F$ is
equipped with a surjective augmentation map $f:F \to \Id$, and
denote by $P:\C \to \C$ the kernel of this map ($P$ is obviously
also an exact functor from $\C$ to iself). The {\em bar-resolution},
or the {\em Godement resolution} associated to the pair $\langle F,f
\rangle$ is constructed as follows. For any $l \geq 1$ we set
\begin{equation}\label{bar.def}
B_l(F,f) = F \circ P^{l-1}:\C \to \C,
\end{equation}
and we define a map $d:B_{l+1}(F,f) \to B_l(F,f)$ as the composition
of the projection $f:F \circ P^l \to \Id \circ P^l = P^l$ and the
embedding $P^l = P \circ P^{l-1} \to F \circ P^{l-1}$. It will be
convenient for us to extend this by setting $B_0=\Id$, $d=f:B_1 = F
\to \Id$, and we will omit $f$ from notation when it is clear from
the context what it is. Then $d^2=0$, and $\langle B_l(F),d \rangle$
is an complex of exact functors from $\C$ to itself. For any object
$E \in \C$, $B_\idot(F)(E)$ is a complex in $\C$ with
$B_0(F)(E)=E$. Moreover, it is easy to check by induction that the
complex $B_\idot(F)(E)$ is actually acyclic, so that $B_{\geq
1}(F)(E)$ is a resolution of the object $E$ -- indeed, by
construction we have a natural embedding of complexes
\begin{equation}\label{boott}
P(B_\idot(F)(E))[1] \to B_\idot(F)(E),
\end{equation}
and the cokernel of this embedding is the acyclic complex $E \to E$.

The bar-construction works without any changes in a more general
situation when $F$ is functor from the category of complexes in $\C$
to itself. In this setting, for any $E \in \C$ the bar complex
$B_\idot(F)(E)$ is naturally a bicomplex; by abuse of notation, we
will denote by $B_\idot(F)(E)$ its total complex.

We note that formally, in order to define the complex $B_\idot(F,f)$
we do not need $P$ to be the kernel of the map $f:F \to \Id$ -- all
we need is a sequence
\begin{equation}\label{PF}
\begin{CD}
0 @>>> P @>>> F @>>> \Id @>>> 0
\end{CD}
\end{equation}
of exact functors which is exact on the left and on the right. In
this yet more general situation, we denote the bar complex by
$B_\idot(F,P)$. This construction is obviously functorial with
respect to sequences of the form \eqref{PF}. If \eqref{PF} is exact
in the middle term, the complex $B_\idot(F,P)=B_\idot(F)$ is
acyclic, if not, then not. However, if $F$ and $P$ in \eqref{PF} are
allowed to be exact functors from the category of complexes in $\C$
to itself, then the condition is weaker: the bar complex
$B_\idot(F,P)(E)$ is acyclic for any $E \in \C$ if the middle
homology of the sequence \eqref{PF} is acyclic -- in other words, if
\eqref{PF} gives an exact triangle after passing to the derived
category.

\subsection{Cyclic categories -- definition and combinatorics.}
The main small categories that we will need in the paper are the
cyclic category $\Lambda$ introduced by A. Connes, and its
generalizations, the so-called paracyclic categories $\Lambda_p$
introduced in \cite{GB}. Let us recall the definitions. Consider the
category $\Cycl$ of linearly ordered sets $M$ equipped with an
order-preserving automorphism $\sigma:M \to M$. Maps in $\Cycl$ are
order-preserving maps which commute with $\sigma$. Consider the set
$\Z$ of all integers as a linearly ordered set with the natural
order, and for any positive $m \in \Z$, denote by $[m]$ the set $\Z$
equipped with an automorphism $\sigma:\Z \to \Z$ given by $a \mapsto a
+ m$. Let $\Lambda_\infty \subset \Cycl$ be the full subcategory
spanned by $[m] \in \Cycl$ for all positive $m \in \Z$. We will
treat $\Lambda_\infty$ is a small category whose set of objects is
the set of positive integers, and for any two positive integers $n,m
\in \Z$, the set
$$
\Lambda_\infty([m],[n])
$$
of maps from $[m]$ to $[n]$ is the set of all order-preserving maps
$f:\Z \to \Z$ such that $f(a+m)=f(a)+n$. In particular, $\sigma$
itself gives an automorphism $\sigma \in \Lambda_\infty([m],[m])$ for
any object $[m]$, and for any $f \in \Lambda_\infty([m],[n])$, we
have $\sigma \circ f = f \circ \sigma$. Therefore setting $f \mapsto f
\circ \sigma$ gives an action of the map $\sigma$ on
$\Lambda_\infty([m],[n])$ which is compatible with compositions, and
we can define a category $\Lambda$ with the same objects as
$\Lambda_\infty$, and with morphism sets given by
$$
\Lambda([m],[n]) = \Lambda_\infty([m],[n])/\sigma.
$$
This is the Connes cyclic category. Moreover, for any positive $p
\in \Z$ we can set
$$
\Lambda_p([m],[n]) = \Lambda_\infty([m],[n])/\sigma^p
$$
and obtain the paracyclic category $\Lambda_p$. For any $p \geq 1$,
we obviously have a natural embedding $i_p:\Lambda_p \to \Lambda$,
$[m] \mapsto [mp]$, and a natural projection $\pi_p:\Lambda_p \to
\Lambda$, $[m] \mapsto [m]$. The projection $\pi_p:\Lambda_p \to
\Lambda$ is a bifibration over $\Lambda$, with all fibers isomorphic
to $\ppt_p$, the groupoid with one object with automorphism group
$\Z/p\Z$. The projection $\rho_p:\Lambda_\infty \to \Lambda_p$ is also
a bifibration, with all fibers isomorphic to $\ppt_\infty$, the
groupoid with one object with automorphism group $\Z$.

The category $\Lambda_\infty$ is self-dual -- the duality functor
$D:\Lambda_\infty \to \Lambda^o_\infty$ is identical on objects, and
for any $f \in \Lambda_\infty([m],[n])$ we set
$$
D(f)(a) = \max\{b \in \Z \mid f(b) \leq a \}, \qquad a \in \Z.
$$
This duality descends to the cyclic category $\Lambda$ and to all
the paracylic categories $\Lambda_p$.

Let $\Delta$ be the simplicial category, that is, the category of
finite linearly ordered sets. Then the opposite category $\Delta^o$
is naturally embedded into $\Lambda_\infty$: for any $[m]$, $[n]$,
the set $\Delta([m],[n])$ is naturally identified with the set of
those $f \in \Lambda_\infty([m],[n])$ that preserve $0 \in \Z$,
$f(0)=0$. Denote this embedding by $j_\infty:\Delta^o \to
\Lambda_\infty$. For every $p \geq 1$, the embedding $j_\infty$
induces an embedding from $\Delta^o$ to $\Lambda_p$, which we denote
by $j_p:\Delta^o \to \Lambda_p$ (simply $j$ if $p=1$). By duality,
we obtain the embeddings $j_p^o:\Delta \to \Lambda_p$. The
bifibration $\pi_p:\Lambda_p \to \Lambda$ splits over $\Delta^o
\subset \Lambda$ by means of the embedding $j_p$ -- in other words,
the embedding $j_p:\Delta^o \to \Lambda_p$ extends to an embedding
$\wt{j}_p:\Delta^o \times \ppt_p \to \Lambda_p$. The composition
$i_p \circ \wt{j}_p:\Delta^o \times \ppt_p \to \Lambda$ factors
through an embedding from $\Delta^o \times \ppt_p$ to $\Delta^o$
which we also denote by $i_p$, by abuse of notation. We can collect
all these data into a commutative diagram
$$
\begin{CD}
\Delta^o @>{j}>> \Lambda\\
@A{i_p}AA  @AA{i_p}A\\
\Delta^o \times \ppt_p @>{\wt{j}_p}>> \Lambda_p\\
@VVV   @VV{\pi_p}V\\
\Delta^o @>{j}>> \Lambda
\end{CD}
$$
of small categories and functors, with Cartesian
squares. Explicitly, the embedding $i_p:\Delta^o \times \ppt_p \to
\Delta^o$ sends $[m] \in \Delta$ to $[p] \times [m]$ with
lexicographical ordering; the automorphism group $\Z/p\Z$ acts by
changing the ordering on the set $[p]$.

As far as the functor categories are concerned, one can also treat
the embedding $j_p:\Delta^o \to \Lambda_p$ as a discrete
cofibration, both for $p \in \Z$ and $p = \infty$. Indeed, let
$\wt{\Delta^o}$ be the total space of the functor $\Lambda_p \to
\Sets$ given by $[n] \mapsto \Lambda_p([1],[n])$. Explicitly,
$\wt{\Delta^o}$ is the category of pairs $\langle [n], a \rangle$ of
the linearly ordered set $\Z$ with the map $\sigma:\Z \to \Z$ and a
fixed element $a \in \Z$ defined modulo $\sigma^p$; maps in
$\wt{\Delta^o}$ are maps in $\Lambda_p$ which preserve the fixed
element. Then the embedding $j_p:\Delta^o \to \Lambda_p$ naturally
factors through an embedding $\Delta^o \to \wt{\Delta^o}$ which
sends $[n]$ to $\langle [n],0\rangle$. It is easy to check that this
embedding $\Delta^o \to \wt{\Delta^o}$ is an equivalence of
categories, and in particular, it induces an equivalence
$\Fun(\wt{\Delta^o},\C) \cong \Fun(\Delta^o,\C)$.

\subsection{Periodicity.}
Fix a Grothendieck abelian category $\C$ and an integer $p \geq
1$, and consider the embedding $j_p:\Delta^o \to \Lambda_p$. 

\begin{lemma}
For any $E \in \Fun(\Lambda_p,\C)$, there exists an isomorphism
$$
j_{p!}j_p^*E \cong E \otimes j_{p!}j_p^*\Z,
$$
this isomorphism is functorial in $E$, and $j_{p!}j_p^*\Z$ is
isomorphic to $\Z_{[1]} \in \Fun(\Lambda_p,\Z\mmod)$.
\end{lemma}

\proof{} Replace $j_p:\Delta^o \to \Lambda_p$ with the discrete
fibration $\xi:\wt{\Delta^o} \to \Lambda_p$; since $\Delta^o$ is
equivalent to $\Delta^o$, we have $j_{p!}j^*_p \cong \xi_!\xi^*$,
and both claims follow from Lemma~\ref{fibr}.
\endproof

By duality, we also obtain a functorial isomorphism
$j^o_{p*}j^{o*}_pE \cong E \otimes j^o_{p*}j^{o*}_p\Z \cong E
\otimes \Z_{[1]}^o$, and since both $\Z_{[1]}$ and $\Z_{[1]}^o$ are
functors into {\em flat} abelain groups, we conclude that both
$j_{p!}j^*_p$ and $j^o_{p*}j^{o*}_p$ are exacts functors from
$\Fun(\Lambda_p,\C)$ to itself. We tautologically have
$H_\idot(\Lambda_p,j_{p!}j^*_pE) \cong H_\idot(\Delta^o,j^*_pE)$.

\begin{lemma}\label{van}
For any $E \in \Fun(\Lambda_p,\C)$, we have
$$
H_\idot(\Lambda_p,j^o_{p*}j^{o*}_pE)=H_\idot(\Lambda_p,E \otimes
\Z_{[1]}^o)=0.
$$
\end{lemma}

\proof{} Since $E \mapsto E \otimes \Z_{[1]}^o$ is an exact functor
in $E$, the homology that we have to study is the derived functor of
the right-exact functor
$$
H_0(\Lambda_p,- \otimes \Z_{[1]}^o):\Fun(\Lambda_p,\C) \to \C.
$$
Therefore it suffices to prove that this right-exact functor
vanishes. Indeed, the constant functor $\Z \in
\Fun(\Lambda_p,\Z\mmod)$ embeds into $\Z_{[1]}^o$ by adjunction, and
the cokernel of this embedding embeds into $\Z_{[2]}^o$; by
adjunction, we have a functorial exact sequence
$$
\begin{CD}
E([2])^{\oplus 2p} @>\rho>> E([1])^{\oplus p} @>>> H_0(\Lambda_p,E
\otimes \Z_{[1]}^o) @>>> 0,
\end{CD}
$$
and it remains to compute the map $\rho$ and to notice that it is
tautologically surjective (in fact, split). We leave it to the
reader.
\endproof

We can now deduce the main homological result on the embedding
$j_p:\Delta^o \to \Lambda_p$.

\begin{lemma}\label{uni.gen}
There exists a map $B_p:\Z_{[1]}^o \to \Z_{[1]}$ such that the
sequence
\begin{equation}\label{seq}
\begin{CD}
0 @>>> \Z @>>> \Z_{[1]}^o @>{B_p}>> \Z_{[1]} @>>> \Z @>>> 0
\end{CD}
\end{equation}
is an exact sequence in $\Fun(\Lambda_p,\Z\mmod)$. The cohomology
algebra $H^\hdot(\Lambda_p,\Z)$ is the free algebra $\Z[u]$ in one
generator $u = u(p)$, and this generator is represented by Yoneda by the
exact sequence \eqref{seq}. The embedding $i_p:\Lambda_p \to
\Lambda$ sends the generator $u(1)$ to $u(p)$.
\end{lemma}

\proof{} The construction of the map $B_p$ for $p=1$ is standard,
see e.g. \cite{FT}. In the case $p > 1$, we notice that
$i_p^*(\Z_{[1]})$ is canonically isomorphic to $\Z_{[1]}$, and
$i_p^*(\Z_{[1]}^o)$ is canonically isomorphic to $\Z_{[1]}^o$. This
gives the map $B_p$ and, by Yoneda, the universal generator $u=u(p)
\in H^2(\Lambda_p,\Z)$. Finally, to show that
$H^\hdot(\Lambda_p,\Z)$ is freely generated by $u$, it suffices to
extend \eqref{seq} to a resolution of the constant object $\Z \in
\Fun(\Lambda_p,\Z\mmod)$ and to notice that by adjunction,
$$
H^\hdot(\Lambda_p,\Z_{[1]}^o) = H^\hdot(\ppt,\Z)=\Z,
$$
while $H^\hdot(\Lambda_p,\Z_{[1]})=0$ by the statement dual to
Lemma~\ref{van}.
\endproof

\begin{corr}\label{per.corr}
For any object $E \in \Fun(\Lambda_p,\C)$, we have a functorial
exact sequence
\begin{equation}\label{per.eq}
\begin{CD}
0 @>>> E @>>> j^o_{p*}j_p^{o*}E @>>> j_{p!}j_p^*E @>>> E @>>>
0
\end{CD}
\end{equation}
and a functorial exact triangle
\begin{equation}\label{per.tr}
\begin{CD}
H_\idot(\Delta^o,j_p^*E) @>>> H_\idot(\Lambda_p,E) @>{u}>>
H_\idot(\Lambda_p,E)[2] @>>>,
\end{CD}
\end{equation}
where $u:H_\idot(\Lambda_p,E) \to H_\idot(\Lambda_p,E)[2]$ is given
by the action of the universal generator $u \in
H^2(\Lambda_p,\Z)$.\endproof
\end{corr}

It will be convenient for us to rewrite \eqref{per.eq} as a
functorial short exact sequence
\begin{equation}\label{per.hash}
\begin{CD}
0 @>>> E[1] @>>> j_p^{\dg}E @>>> E @>>> 0
\end{CD}
\end{equation}
in the category of complexes of objects in $\Fun(\Lambda_p,\C)$,
where $E[1]$ is, as usually, the complex consisting of $E$ placed in
degree $-1$, and $j_p^{\dg}E$ is the complex $j^o_{p*}j_p^{o*}E \to
j_{p!}j_p^*E$ placed in degrees $0$ and $-1$. Then \eqref{per.tr}
amounts to a canonical isomorphism $H_\idot(\Delta^o,j_p^*E) \cong
H_\idot(\Lambda_p,j_p^{\dg}E)$ and the long exact sequence associated
to \eqref{per.hash}. We note that $j_p^{\dg}$ is an exact functor
from $\Fun(\Lambda_p,\C)$ to the category of complexes in
$\Fun(\Lambda_p,\C)$; in particular, it trivially extends to a
functor $j_p^{\dg}:\D^-(\Lambda_p,\C) \to \D^-(\Lambda_p,\C)$ on the
derived categories. If $p=1$, one usually denotes
\begin{align}\label{hh.hc}
\begin{split}
H_\idot(\Delta^o,j^*E) &= HH_\idot(E),\\
H_\idot(\Lambda,E) &=HC_\idot(E),
\end{split}
\end{align}
and calls them the {\em Hochschild} and {\em cyclic} homology of the
cyclic object $E$. The map $u$ is the Connes' periodicity map, and
the exact triangle \eqref{per.tr} is known as the Connes' exact
sequence.

\subsection{Compatibility for paracyclic embeddings.}\label{para.sub}
We finish this section with a somewhat surprising result which shows
that for any cyclic object $E$ and any integer $p$, one can find a
$p$-cyclic object with the same homology (this is the key point in
our construction of the non-commutative Cartier operator).

\begin{prop}\label{i.p}
For any Grothendieck abelian category $\C$ and any object $E \in
\Fun(\Lambda,\C)$, the natural map
$$
i_{p!}:H_\idot(\Lambda_p,i_p^*E) \to H_\idot(\Lambda,E)
$$
is a quasiisomorphism.
\end{prop}

\proof{} By Lemma~\ref{uni.gen} the map $i_{p!}$ is compatible with
the periodicity map $u$; therefore by \eqref{per.tr} it suffices to
prove that the natural map
$$
i_{p!}:H_\idot(\Delta^o,i_p^*E') \to H_\idot(\Delta^o,E')
$$ 
is a quasiisomorphism, where $E' = j^*E$, and $i_p$ is now treated
as the embedding from $\Delta^o$ to itself. We will in fact prove
this claim for any $E' \in \Fun(\Delta^o,E')$. To do this, it
suffices to consider the case when $E'$ goes over a family of
generators of the category $\Fun(\Delta^o,\C)$; to construct
such a family, it suffices to take a generator $S$ of the category
$\Fun(\Delta^o,\Z\mmod)$ and tensor it with all objects in $\C$
in turn. Thus we are reduced to the case $\C = \Z\mmod$, $E' = S$ is
a generator of $\Fun(\Delta^o,\Z\mmod)$. 

Consider the category of bimodules over the polynomial algebra
$B=\Z[t]$ in one variable $t$ over $\Z$, and let $S$ be the standard
bar-resolution of the diagonal bimodule $B$ (we have
$S([n])=B^{\otimes (n+1)}$). This is a simplicial abelian group
which in addition inherits the standard grading from the polynomial
algebra $B$. It is easy to check that the degree-$l$ component $S^l$
of the simplicial $\Z$-module $S$ is isomorphic to the functor
$\Z_{[l]}$ (the standard $l$-simplex); therefore $S = \bigoplus S^l$
is indeed a generator of the category $\Fun(\Delta^o,\Z\mmod)$. On
the other hand, it is well-known that the homology of $\Delta^o$
with coefficients in a simplicial abelian group can be computed by
the standard simplicial complex, and since $S$ is a resolution, we
have
$$
H_0(\Delta^o,S) = B
$$
and $H_i(\Delta^o,S)=0$ for $i \geq 1$. It remains to notice
that $S$ is a simplicial flat $B$-bimodule, and
$$
i_p^*S \cong S \otimes_B S \otimes_B \dots
\otimes_B S,
$$
where the product contains $p$ terms. Since $B \otimes_B B \cong B$,
we conclude that $i_p^*S$ is another resolution for the same
bimodule $B$, and the natural map $H_\idot(\Delta^o,i_p^*S) \to
H_\idot(\Delta^o,S)$ is indeed a quasiisomorphism, as required.
\endproof

\section{Filtrations on cyclic homology.}\label{filtr}

\subsection{The Hodge filtration.} Fix a Grothendieck abelian
category $\C$, and consider an arbitrary object $E \in
\Fun(\Lambda_p,\C)$. Iterating \eqref{per.eq}, one obtains a
canonical resolution $E_\idot$ of the object $E$; using this
resolution, one can refine \eqref{per.hash} in the following
way. Consider the stupid filtration on the complex $E_\idot$, and
let $F^lE_\idot$, $l\geq 0$, denote the $2l$-th term of this stupid
filtration. Then $\langle E,F^\hdot \rangle$ is a filtered complex,
and it defines an object in the filtered derived category
$\DF^-(\Lambda_p,\C)$. The following is a reformulation of
Corollary~\ref{per.corr}.

\begin{lemma}
The correspondence $E \mapsto \langle E_\idot,F^\hdot \rangle$
extends to a functor 
$$
\D^-(\Lambda_p,\C) \to \DF^-(\Lambda_p,\C),
$$
and we have a functorial exact triangle
$$
\begin{CD}
j_p^{\dg}E @>>> E @>{u}>> E[2] @>>>,
\end{CD}
$$
where the filtration in the right-hand side is shifted by $1$, and
$u$ is the canonical periodicity map.\endproof
\end{lemma}

\begin{defn}
The filtration $F^\hdot$ on the object $E \in \D^-(\Lambda_p,\C)$ is
called the {\em Hodge filtration}.
\end{defn}

By the standard formalism of filtered complexes, the Hodge
filtration induces a spectral sequence in homology which starts with
$HH_\idot(E)[v]$, the space of polynomials in one variable $v =
u^{-1}$ of degree $-2$, and converges to $HC_\idot(E)$. We call it
the {\em Hodge-to-de Rham}, or {\em Hochschild-to-cyclic} spectral
sequence. It is this spectral sequence whose degeneration we are
going to study in the rest of the paper.

Assume now that the integer $p \geq 1$ is actually an odd prime, and
assume that the category $\C$ is $k$-linear over a field $k$ of
characteristic $p$. In this case, there is another way of looking at
the Hodge filtration. Namely, consider the projection
$\pi=\pi_p:\Lambda_p \to \Lambda$ (from now on, we will fix $p$ and
drop it from the notation whenever there is no danger of
confusion). By definition, for any $E \in \Fun(\Lambda_p,\C)$ we
have
$$
H_\idot(\Lambda_p,E) \cong H_\idot(\Lambda,\Ll^\hdot\pi_!(E)),
$$
where $\Ll^\hdot\pi_!:\D^-(\Lambda_p,\C)\to\D^-(\Lambda,\C)$ is the
derived functor of the direct image functor $\pi_!$. We also have
$H^\hdot(\Lambda,\Rr^\hdot\pi_*k) \cong H^\hdot(\Lambda_p,k)$;
consider the associated spectral sequence (the Hochschild-Serre
spectral sequence for the fibration $\pi:\Lambda_p \to
\Lambda$). Recall that the fibers of the fibration $\pi:\Lambda_p
\to \Lambda$ are all isomorphic to the groupoid $\ppt_p$, so that
taking $\Ll^\hdot\pi_!$ amounts to taking homology of the finite
group $\Z/p\Z$. Recall also that the cohomology algebra
$H^\hdot(\Z/p\Z,k)$ of the cyclic group $\Z/p\Z$ with coefficients
in the trivial module $k$ is the free graded-commutative algebra
generated by a generator $\eps$ of degree $1$ and a generator $u$ of
degree $2$.

\begin{lemma}\label{HS}
The graded cyclic $k$-vector space $\Rr^\hdot\pi_*k$ is a free
module over the cohomology algebra
$H^\hdot(\Z/p\Z,k)=k[u]\langle\eps\rangle$; in particular, for any
$l \geq 0$ we have $\Rr^l\pi_*k \cong k \in
\Fun(\Lambda,k\mmod)$. The Hochschild-Serre spectral sequence for
the fibration $\pi:\Lambda_p \to \Lambda$ degenerates at $E_3$, we
have
$$
E_2^{r,s} = H^s(\Lambda,\Rr^rk) =
\begin{cases}
k, \qquad s \text{ is even},\\
0, \qquad s \text{ is odd},
\end{cases}
$$
and the differential $d_2:E_2^{r,s} \to E_2^{r-1,s+2}$ vanishes for
even $r$ and induces an isomorphism for odd $r$.
\end{lemma}

\proof{} To prove the first claim, it suffices to evaluate
$\Rr^\hdot\pi_*k([n])$ for all $[n] \in \Lambda$ and to notice that by
base change, we have 
$$
\Rr^\hdot\pi_*k([n]) = H^\hdot(\Z/p\Z,k([n]))=H^\hdot(\Z/pZ,k).
$$
To compute $E_2^{r,s}$, combine this with Lemma~\ref{uni.gen}. Again
by Lemma~\ref{uni.gen}, the class $\eps \in E_2^{1,0}$ cannot
survive in the $E_\infty$-term of the spectral sequence; therefore
$d_2(\eps)$ generates $E_2^{0,2}$. Once again by
Lemma~\ref{uni.gen}, the class $u \in E_2^{2,0}$ must survive in
$E_\infty$, so that in particular $d_2(u)=0$. The claim about $d_2$
follows by multiplicativity. This implies that $E_3^{r,s}=0$ unless
$s=0,1$ and $r=2l$ is even, and the spectral sequence degenerates at
$E_3$ by dimension reasons.
\endproof

\begin{remark}\label{nonspl}
This Lemma shows, in particular, that $\Lambda_p$ considered as a
gerbe over $\Lambda$ is {\em not} split -- conversely, its class in
$H^2(\Lambda,\Z/p\Z)$ is the generator class $u$.
\end{remark}

\begin{lemma}\label{acy}
For any $E \in \Fun(\Lambda_p,\C)$, the complex $j_p^{\dg}(E) \in
\Fun(\Lambda_p,\C)$ consists of objects acyclic for the functors
$\pi_!,\pi_*:\Fun(\Lambda_p,\C) \to \Fun(\Lambda,\C)$.
\end{lemma}

\proof{} By the projection formula, it suffices to prove that
$j_{p!}(k),j_{p*}(k) \in \Fun(\Lambda_p,k\mmod)$ are acyclic for
$\pi_!$ and $\pi_*$. By Proposition~\ref{fibr}, this has to be
checked for every fiber $\pi_{[n]}$ of the bifibration
$\pi:\Lambda_p \to \Lambda$. This fiber is equivalent to $\ppt_p$,
so that we have to check that $j_{p!}(k)([n])$ and $j_{p*}(k)([n])$
are regular representations of the group $\Z/p\Z$. This is
immediate.
\endproof

By virtue of Lemma~\ref{acy}, one can compute $\Ll^\hdot\pi_!E$ by
using the standard periodic resolution $E_\idot$ assembled out of
the complexes $j_p^{\dg}(E)$. This defines a filtration $F^\hdot$ on
$\Ll^\hdot\pi_!E$ and equips it with the structure of a filtered
complex in $\Fun(\Lambda,\C)$; by abuse of terminology, we will also
call the filtration $F^\hdot\Ll^\hdot\pi_!E$ the {\em Hodge
filtration}.

\subsection{The conjugate filtration.}\label{tate.sub}
It turns out, however, that the complex $\Ll^\hdot\pi_!E$ carries
another canonical filtration, which goes in the opposite direction;
we call it the conjugate filtration. We start with the following.

\begin{lemma}\label{tate}
An object $F \in \Fun(\Lambda_p,k\mmod)$ is acyclic for the functor
$\pi_!$ if and only if it is acyclic for the functor $\pi_*$;
moreover, there exists a functorial map $T:\pi_!F \to \pi_*F$ which
is an isomorphism for acyclic $F$. For any $k$-linear Grothendieck
category $\C$, any $E \Fun(\Lambda_p,\C)$ and any acyclic $F \in
\Fun(\Lambda_p,k\mmod)$, $F \otimes E \in \Fun(\Lambda_p,\C)$ is
acyclic for both $\pi_!$ and $\pi_*$, and $T:\pi_!(E \otimes F) \to
\pi_*(E \otimes F)$ is an isomorphism.
\end{lemma}

\proof{} As in the proof of Lemma~\ref{acy}, it suffices by
Proposition~\ref{fibr} to check the statement for every fiber
$\pi_{[n]}$ of the bifibration $\pi:\Lambda_p \to \Lambda$; since
$\pi_{[n]} \cong \ppt_p$, this reduces to standard facts about the
homology of finite groups -- in this case, the group is
$\Z/p\Z$. The map $T$ is the Tate trace map: for every $[n] \in
\Lambda$, $T([n]) = 1 + \sigma + \dots + \sigma^{p-1}$ is the
averaging over $\Z/p\Z$, where $\sigma \in \Z/p\Z$ is the generator.
\endproof

Therefore we can take a left acyclic resolution $I_\idot$ and a
right acyclic resolution $I^\hdot$ of the constant functor $k \in
\Fun(\Lambda_p,k\mmod)$ and combine them into an unbouded complex
$I_{\#}$ by taking the cone of the natural map $I_\idot \to k
I^\hdot$ (we normalize the grading so that $I_{\#}$ is an extension
of $I_\idot[1]$ by $I^\hdot$).

\begin{defn}
For $k$-linear Grothendieck category $\C$ and for any $E \in
\Fun(\Lambda_p,\C)$, the unbounded complex
$$
\pi_{\#}(E) = \pi_*(E \otimes I_{\#}) \in D(\Lambda,\C)
$$
is called the {\em relative Tate homology complex} of $E$ with
respect to $\pi:\Lambda_p \to \Lambda$.
\end{defn}

The Tate homology complex $\pi_{\#}(E) \in D(\Lambda,\C)$ is
obviously independent of the choice of resolutions $I_\idot$,
$I^\hdot$ and functorial in $E$. We have a functorial exact triangle
$$
\begin{CD}
\Ll^\hdot\pi_!(E)[1] @>>> \pi_{\#}(E) @>>> \Rr^\hdot\pi_*(E)
@>>>
\end{CD}
$$
in $\D(\Lambda,k\mmod)$ and a periodicity map $u:\pi_{\#}(E) \to
\pi_{\#}(E)[2]$ which is compatible with the periodicity maps on
$\Ll^\hdot\pi_!(E)$ and $\Rr^\hdot\pi_*(E)$.

\begin{defn}
The {\em conjugate filtration} $W_\idot$ on the complex
$\pi_{\#}(E)$ is defined as its canonical truncation
$$
W_l\pi_{\#}(E) = \tau_{-2l-1}\pi_{\#}(E).
$$
\end{defn}

This is also functorial and independent of choices. Moreover, by
construction we see that $W_l$ for $l \geq 1$ is actually a
filtration on $\Ll^\hdot\pi_!(E)[1] \cong \pi_!(E \otimes
I_\idot)[1]$. We will adopt this point of view to avoid dealing with
unbounded complexes. However, for various reasons it will be more
convenient to consider also the canonical truncation $\tau_{\leq
0}\pi_{\#}(E)$ equipped with the filtration induced by $W_\idot$,
and we will denote this filtered complex by
$$
\pi_{\flat}(E) = \tau_{\leq 0}\pi_{\#}(E).
$$
The complex $\pi_{\flat}(E)$ is an extension of
$\Ll^\hdot\pi_!(E)[1]$ by $\pi_*(E)$.

We note that if we compute the complex $\pi_{\#}(E)$ by using the
standard periodic resolution $E_{\#}$, then the periodicity map $u$
is induced by an invertible map of complexes $u:E_{\#} \to
E_{\#}[2]$. We denote the inverse map by
$$
v = u^{-1}:E_{\#}[2] \to E_{\#}.
$$
This map shifts the conjugate filtration by $1$ -- we actually have
a series of maps
\begin{equation}\label{perio}
v: W_lE_{\#}[2] \to W_{l+1}E_{\#}.
\end{equation}

\subsection{Periodicity in the bar resolutions.}\label{bar.sub}
To analyse the conjugate filtration, we will need to use several
acyclic resolution in Lemma~\ref{tate}. We use the machinery of bar
resolutions described in Subsection~\ref{prelim}.  Consider the
category $\Fun(\Lambda_p,\C)$ of $p$-cyclic objects in a
Grothendieck abelian category $\C$. There are two natural exact
functors from $\Fun(\Lambda_p,\C)$ to itself that can play the role
of $F$ in \eqref{bar.def}: firstly, we can take the functor
$j_{p!}j_p^*$, secondly, we can take the functor $j_p^{\dg}$ (the
second one takes values in complexes in
$\Fun(\Lambda_p,\C)$). Moreover, in the second case we can take $P =
\Id[1]$, with the natural embedding $\Id[1] \to j_{p*}j^*_p[1] \to
j_p^{\dg}$. This gives three possible functorial resolution related
by the following natural maps
\begin{equation}\label{3.way}
\begin{CD}
B_\idot(j_{p!}j^*_p)(E) @>>> B_\idot(j_p^{\dg})(E) @<<<
B_\idot(j_p^{\dg},\Id[1])(E).
\end{CD}
\end{equation}
All these maps are quasiisomorphisms, and all the complexes are
acyclic. To simplify notation, denote
$B_\idot(j_{p!}j^*_p)=B_\idot^!$, $B_\idot(j_p^{\dg}) =
B_\idot^{\dg}$ and $B_\idot(j_p^{\dg},\Id[1])=\bB^{\dg}_\idot$. By
Lemma~\ref{acy}, for any $l \geq 1$, all three objects $B_l^!(E)$,
$B_l^{\dg}(E)$ and $\bB_l^{\dg}(E)$ are acyclic for the functor
$\pi_*$, so that all three resolutions can be used in
Lemma~\ref{tate} to compute the relative Tate homology complex
$\pi_{\flat}(E)$. For any $E \in \Fun(\Lambda_p,\C)$,
\begin{equation}\label{tau.0}
\pi_*(B_\idot^!(E)) \cong \pi_*(B_\idot^{\dg}(E)) \cong
\pi_*(\bB_\idot^{\dg}(E)) \cong \pi_{\flat}(E).
\end{equation}
The complex $\bB_\idot^{\dg}(E)$ is nothing but the standard
periodic resolution of $E$ composed of the complexes
$j_p^{\dg}(E)[2l]$, $l \geq 0$, and the filtration on
$\bB_\idot^{\dg}(E)$ is the conjugate filtration. The canonical map
\eqref{boott} for this complex is a map
\begin{equation}\label{bB}
v:\bB_\idot^{\dg}(E)[2] \to \bB_\idot^{\dg}(E),
\end{equation}
which coincides with the periodicity map \eqref{perio}. For the
complex $B_\idot^{\dg}(E)$, the map \eqref{boott} is the map
\begin{equation}\label{B.dg}
\overline{j}_p^{\dg}(B_\idot^{\dg}(E))[1] \to B_\idot^{\dg}(E),
\end{equation}
where we denote by $\overline{j}_p^{\dg}$ the kernel of the natural
map $j_p^{\dg} \to \Id$. Since the natural map $\Id[1] \to
\overline{j}^{\dg}_p$ is a quasiisomorphism, the natural map
$$
B^{\dg}_\idot(E) \to \bB^{\dg}_\idot(E)
$$
is a quasiisomorphism for every $E \in \Fun(\Lambda_p,\C)$. To see
the periodicity map $v$, we have to compose the map \eqref{B.dg}
with the natural quasiisomorphism $B^{\dg}_\idot[1](E) \to
\overline{j}_p^{\dg}(B_\idot^{\dg}(E))$ induced by the
quasiisomorphism $\Id[1] \to \overline{j}^{\dg}_p$.

In Section~\ref{horror}, we will need to use the complex
$B^!_\idot(E)$, and in particular, to interpret the periodicity in
terms of this complex. To do this, we note that $\pi_* \circ j_{p*}
\cong j_*$, and by Lemma~\ref{acy} and Lemma~\ref{tate} we have
$\pi_* \circ j_{p!}  \cong \pi_! \circ j_{p!} \cong j_!$. Therefore
$\pi_*(j_p^{\dg}k) \cong j^{\dg}k$, and by adjunction, we have a map
\begin{equation}\label{phi}
\phi:\pi^*(j^{\dg}k) \to j^{\dg}_pk.
\end{equation}
Both sides are complexes of length $2$ with homology in degree $0$
and $1$ equal to the constant functor $k$; the right-hand side is a
Yoneda representation of the periodicity class $u(p) \in
H^2(\Lambda_p,k)$, and the left-hand side represents $\pi^*(u) \in
H^2(\Lambda_p,k)$, which is equal to $0$ by Lemma~\ref{HS}.  The map
$\phi$ is an isomorphism on homology in degree $1$, and since
$\pi^*(u)=0$, it is trivial on homology in degree $0$. Thus the map
$\phi$ actually maps $\pi^*(j^{\dg}k)$ into $\overline{j}^{\dg}_pk
\subset j^{\dg}_pk$ and thus extends the canonical embedding $k[1]
\to \overline{j}_p^{\dg}k$. Composing $\phi$ with the canonical map
\eqref{B.dg}, we extend the periodicity map \eqref{bB} to a map
\begin{equation}\label{phi.B}
v:\pi^*(j^{\dg}(k)) \otimes B^{\dg}_\idot(E)[1] \to B^{\dg}_\idot(E).
\end{equation}
We will now denote by $k(1) \subset j_!k \in \Fun(\Lambda_p,k\mmod)$
the kernel of the canonical map $j_!k \to k$, and for every $E \in
\Fun(\Lambda,\C)$, or $E \in \Fun(\Lambda_p,\C)$, we will denote
$E(1) = E \otimes k(1)$, resp. $E(1) = E \otimes \pi^*k(1)$. The map
$\phi$ in \eqref{phi} induces in particular a map $\phi:E(1) \to
j_{p!}j_p^*E$, which actually goes into the kernel of the canonical
map $j_{p!}j_p^*E \to E$. Composing this map $\phi$ with the
canonical map \eqref{boott} for the complex $B_\idot^!(E)$, we
obtain a functorial map
\begin{equation}\label{bV}
\overline{v}:B_\idot^!(E)(1)[1] \to B_\idot^!(E).
\end{equation}
Now, the object $k(1) \in \Fun(\Lambda,k\mmod)$ is also the cokernel
of the canonical map $k \to j_*k$. Therefore by Lemma~\ref{van} we
have
$$
H_\idot(\Lambda,E(1)) \cong H_\idot(\Lambda,E)[1]
$$
for any $E \in \Fun(\Lambda,\C)$. More generally, by the projection
formula, we have
\begin{equation}\label{1=1}
H_\idot(\Lambda,\pi_*(E(1))) \cong
H_\idot(\Lambda,\pi_*(E)(1)) \cong H_\idot(\Lambda,\pi_*E)[1]
\end{equation}
for any $E \in \Fun(\Lambda_p,\C)$; 

\begin{lemma}\label{bar.perio}
For any $E \in \Fun(\Lambda_p,\C)$, the map
$$
H_\idot(\Lambda,\pi_*B_\idot^!(E)(1)[1]) \to
H_\idot(\Lambda,\pi_*B_\idot^!(E))
$$
induced by the map $\overline{v}$ from \eqref{bV} and the map
$$
H_\idot(\Lambda,\pi_*\bB_\idot^{\dg}(E)[2]) \to
H_\idot(\Lambda,\pi_*\bB_\idot^{\dg}(E))
$$
induced by the periodicity map $v$ from \eqref{bB} become equal
under the identification \eqref{1=1}.
\end{lemma}

\proof{} By construction, the diagram \eqref{3.way} extends to a
diagram
\begin{equation}\label{3.way.2}
\begin{CD}
B_\idot^!(E)[1](1) @>>> B_\idot^{\dg}[1](E) \otimes
\pi^*\overline{j}^{\dg}k @<<< \bB_\idot^{\dg}(E)[2]\\
@V{\overline{v}}VV @VV{v}V @VV{v}V\\
B_\idot^!(E) @>>> B_\idot^{\dg}(E) @<<< \bB_\idot^{\dg}(E),
\end{CD}
\end{equation}
and by the projection formula, it suffices to prove that for any
complex $E'$ in $\Fun(\Lambda,\C)$ -- in particular, for $E' =
\pi_*B_\idot^!(E)$ and such -- the natural maps
$$
\begin{CD}
H_\idot(\Lambda,E'(1)) @>>> H_\idot(\Lambda,E' \otimes
\overline{j}^{\dg}(k)) @<<< H_\idot(\Lambda,E'[1])
\end{CD}
$$
are quasiisomorphisms. Indeed, the map on the right-side is a
quasiisomorphism already in $\Fun(\Lambda,\C)$, while the cone of
the map on the left-hand side is $H_\idot(\Lambda,j_*j^*E')$, which
is trivial by Lemma~\ref{van}.
\endproof

As an application, Lemma~\ref{bar.perio} allows to see the conjugate
filtration in the bar resolutions $B_\idot^!(E)$ and
$B_\idot^\dg(E)$. Namely, by definition $\bB_\idot^\dg(E)$ carries
the conjugate filtration $W_\idot$, and in particular, we have
$W_1\bB_\idot^\dg(E) = \tau_{\leq -1}\bB^\dg_\idot(E)$. Set
$W_1B_\idot^!(E) = \tau_{\leq -1}\B_\idot^!(E)$,
$W_1B_\idot^\dg(E) = \tau_{\leq -1}\B_\idot^\dg(E)$, and note that
all the vertical maps in \eqref{3.way.2} are injective. Therefore we
can define inductively
\begin{align*}
W_lB_\idot^!(E) &= v\left(W_{l-1}B_\idot^!(E)[1](1)\right),\\
W_lB_\idot^\dg(E) &= v\left(W_{l-1}B_\idot^\dg(E)[1](1) \otimes
\pi^*\overline{j}^{\dg}k\right).
\end{align*}
This turns $B_\idot^!(E)$ and $B_\idot^\dg(E)$ into filtered
complexes, and all the maps in \eqref{3.way} are filtered maps.

\begin{corr}\label{bar.perio.corr}
The natural map $\pi_*(\bB_\idot^\dg(E)) \to \pi_*(B_\idot^\dg(E))$
is a filtered quasiisomorphism, while the natural map
$\pi_*(B_\idot^!(E)) \to \pi_*(B_\idot^\dg(E))$ becomes a filtered
quasiisomorphism aftr applying the cyclic homology functor
$HC_\idot$.
\end{corr}

\proof{} By Lemma~\ref{bar.perio}, it suffices to notice that the
map $k[1] \to j^\dg k$ is a quasiisomorphism in
$\Fun(\Lambda,k\mmod)$, while the map $k(1) \to j^\dg k$ becomes a
filtered quasiisomorphism after applying $HC_\idot(-) =
H_\idot(\Lambda,-)$.
\endproof

\subsection{Tight $\Z/p\Z$-modules.}\label{tight.sub}
We will now use the bar resolutions to compute, under some
assumptions on the $p$-cyclic object $E \in \Fun(\Lambda_p,\C)$, the
associated graded quotients of the conjugate filtration on
$\pi_{\flat}(E)$.

Let $V$ be a vector space over the field $k$ -- or, more generally,
an object in a $k$-linear Grothendieck category $\C$ -- equipped
with an action of the group $\Z/p\Z$. Denote by $\kk{V}$ the complex
\begin{equation}\label{Vdag}
\begin{CD}
V_{\Z/p\Z} @>{T}>> V^{\Z/p\Z}
\end{CD}
\end{equation}
placed in degrees $0$ and $-1$ (here $T$ is the Tate trace map from
coinvariants to invariants with respect to the group $\Z/p\Z$). The
natural map $V^{\Z/p\Z} \to V \to V_{\Z/p\Z}$ factors through a map
$$
\phi:H_0(\kk{V}) \to H_1(\kk{V}).
$$
The complex $\kk{V}$ is actually a piece of the standard periodic
complex which computes the homology $H_\idot(\Z/p\Z,V)$, so that we
have
$$
\begin{aligned}
H_0(\kk{V}) &\cong H_{odd}(\Z/p\Z,V),\\
H_1(\kk{V}) &\cong H_{even}(\Z/p\Z,V),
\end{aligned}
$$
and the map $\phi$ is given by the action of the standard generator
$\eps$ of the first cohomology group $H^1(\Z/p\Z,k)$.

\begin{defn}\label{tight}
The representation $V$ of the group $\Z/p\Z$ is called {\em tight}
if the map $\phi$ is an isomorphism.
\end{defn}

For example, a trivial representation $k$ is tight, and so is the
regular representation $k[\Z/p\Z]$. In other words, a representation
is tight if $\eps:H_{odd}(\Z/p\Z,V) \to H_{even}(\Z/p\Z,V)$ is an
isomorphism. We note that since $\eps^2=0$, this automatically
implies that $\eps:H_{even}(\Z/p\Z,V) \to H_{odd}(\Z/p\Z,V)$ is a
trivial map. For any tight $\Z/p\Z$-module $V$ in a category $\C$,
we will denote by $\ii{V}$ the object $H_0(\kk{V}) \cong
H_1(\kk{V}) \in \C$.

Assume given an object $E \in \Fun(\Lambda_p,\C)$; then $E([m])$ is
a $\Z/p\Z$-module in $\C$ for every $[m] \in \Lambda_p$, and all the
complexes $\kk{E([m])}$ fit together into a single complex
$$
\pi_!E \to \pi_*E,
$$
which we denote by $\kk{E} \in \Fun(\Lambda,\C)$. The maps $\phi$
for various $E([m])$ fit into a single map
\begin{equation}\label{phi.kk}
\phi = \phi_E:H_0(\kk{E}) \to H_1(\kk{E}).
\end{equation}

\begin{defn}\label{tight.cat}
An object $E \in \Fun(\Lambda_p,\C)$ is called {\em tight} if the
map $\phi_E$ is an isomorphism -- or, equivalently, if for any $[m]
\in \Lambda$ the object $E([m])$ is tight with respect to $\ppt_p
\times [m] = \pi_{[m]} \subset \Lambda_p$. For a tight object $E \in
\Fun(\Lambda_p,\C)$, we denote by $\ii{E} \in \Fun(\Lambda,\C)$ the
object $H_0(\kk{E}) \cong H_1(\kk{E})$.
\end{defn}

\begin{lemma}\label{adag}
For any tight $E \in \Fun(\Lambda_p,\C)$, the periodicity map $v$
induces a filtered exact triangle
\begin{equation}\label{tri.1}
\begin{CD}
\pi_{\flat}(E)[1] @>{v}>> \pi_{\flat}(E)
@>>> \kk{E} @>>>
\end{CD}
\end{equation}
in $\DF^-(\Lambda,\C)$, where the conjugate filtration on the left
is shifted by $1$, and the complex $\kk{E}$ is equipped with the
canonical filtration. Moreover, for any $l \geq 1$ the map $v$
induces a filtered exact triangle
\begin{equation}\label{tri.2}
\begin{CD}
W_{l+1}\pi_{\flat}E @>>> W_l\pi_{\flat}E @>>> j^{\dg}\ii{E}[2l] @>>>,
\end{CD}
\end{equation}
where again the conjugate filtration on the left is shifted by $1$,
and the complex $j^{\dg}\ii{E}$ is equipped with a trivial one-step
filtration of degree $l$.
\end{lemma}

\proof{} Recall that we have filtered quasiisomorphisms
$$
\pi_{\flat}(E) \cong \pi_*(\bB^{\dg}_\idot(E)) \cong
\pi_*(B^{\dg}_\idot(E)).
$$
By definition, $\pi_*(\bB^{\dg}_0(E)) \cong \pi_*(E)$, and
$H_1(\pi_*(\bB^{\dg}_{\geq 1}(E)) \cong \pi_!(E)$, so that we have a
natural projection $\pi_*(\bB^{\dg}_\idot(E)) \to \kk{E}$, and it
obviously fits into a sequence
$$
\begin{CD}
0 @>>> \pi_*(\bB^{\dg}_\idot(E))[2] @>{v}>>
\pi_*(\bB^{\dg}_\idot(E)) @>>> \kk{E} @>>> 0
\end{CD}
$$
of filtered complexes. This sequence is not exact in the middle
term, but its homology is an acyclic complex with trivial one-step
filtration, so that after passing to the filtered derived category
we obtain an exact triangle. This is the triangle \eqref{tri.1}. To
prove \eqref{tri.2}, we may assume by induction that $l=1$. The
periodicity map $v$ is filtered by construction and identical in
degrees $\neq 1$; we have to prove that $\gr^W_1\pi_{\flat}(E) \cong
j^{\dg}\ii{E}$. Indeed, by the projection formula, the canonical
filtered map \eqref{phi.B} induces a map
$$
j^{\dg}k \otimes \pi_*(B^{\dg}_\idot(E)) \cong
j^{\dg}\pi_*(B^{\dg}_\idot(E)) \to
\pi_*(B^{\dg}_\idot(E)),
$$
which in particular gives a map
$$
\phi:j^{\dg}(\gr^W_0\pi_*(B^{\dg}_\idot(E))) \to
\gr^W_1\pi_*(B^{\dg}_\idot(E)).
$$
Since $\gr^W_0\pi_*(B^{\dg}_\idot(E)) \cong H_0(\kk{E})$ by
\eqref{tri.1}, it suffices to prove that this map $\phi$ is a
quasiisomorphism. Both sides are complexes with non-trivial homology
in degrees $0$ and $1$ only, and, again by \eqref{tri.1}, the
homology in degree $1$ on both sides is equal to $H_0(\kk{E})$ and
identified by $\phi$. In degree $0$, the homology on the right-hand
side is $H_1(\kk{E})$, and the homology on the left-hand side is
again $H_0(\kk{E})$. It is elementary to check that the induced map
$\phi:H_0(\kk{E}) \to H_1(\kk{E})$ is the same map as in
\eqref{phi.kk}; since $E$ is assumed tight, $\phi$ is an
isomorphism.
\endproof

As a corollary, we see that for any tight $E \in
\Fun(\Lambda_p,\C)$, the associated graded quotients
$\gr^W_lHC_\idot(E)$ are isomorphic to the cyclic homology
$HC_\idot(j^{\dg}\ii{E})$, which is by definition equal to the {\em
Hochschild} homology $HH_\idot(\ii{E})$. This is the
linear-algebraic origin of our non-commutative Cartier isomorphism.

\section{Conjugate filtration in detail.}\label{horror}

We will now investigate the conjugate filtration on cyclic homology
in some detail. Unfortunately, it seems that apart from
Lemma~\ref{adag}, nothing can be said for general $p$-cyclic
objects, even for those which are tight. Thus in general, we can
identify the associated graded quotients of the conjugate filtration
but cannot study the extension data between these quotients. In
order to have any control over the extension data, we need to impose
an additional symmetry on our objects, which we will call {\em
polycyclic}.

\subsection{Polycyclic groups.}\label{poly.grp}
For any finite set $S$, denote by $G_S$ the abelian group
$\Z/p\Z[S]$ -- the free $\Z/p\Z$-module generated by $S$. Fix a
finite set $S$ and assume given a $G_S$-module in a $k$-linear
abelian category $\C$ over a field $k$ of characteristic $p$. For
every subset $I \subset S$, denote
$$
V_I = H_0(G_I,H^0(G_{\overline{I}},V)),
$$
where $\overline{I} = S \setminus I$ is the complement to $I \subset
S$. For any two subsets $I \subset I' \subset S$, we have natural
transition maps
$$
\tau_{I,I'}:V_{I'} \to V_I \qquad \sigma_{I,I'}:V_I \to V_{I'}
$$
induced by the natural embedding $H^0(G_{\overline{I}},V) \to
H^0(G_{\overline{I}'},V)$ and the trace map
$H^0(G_{\overline{I'}},V) \to H^0(G_{\overline{I}},V)$. Fix a linear
order on $S$, so that $S = \{i, 1 \leq i \leq n\}$, where $n=|S|$ is
the number of elements in $S$. For any $I \subset S$, $i \in
\overline{I}$, let $l(i,I)$ be the number of elements in $I$ which
are less than $i$. Then we define a map
$$
d_I = \sum_{i \in \overline{I}}(-1)^{l(i,I)}\tau_{I,I \cup \{i\}}:V_I \to
\bigoplus_{i \in \overline{I}}V_{I \cup \{i\}},
$$
and if we set
$$
\kkk{V}{l} = \bigoplus_{|I|=l}V_I, \qquad d_l = \bigoplus_{|I|=l}d_I,
$$
then $d_{l+1} \circ d_l=0$, so that $\langle \kk{V},d_\idot
\rangle$ becomes a complex of length $|S|+1$. For any $i \in S$,
define a map $v_i:\kk{V} \to \kkk{V}{\idot+1}$ by
\begin{equation}\label{phi.i}
\phi_i = \sum_{i \in I} (-1)^{l(i,\overline{I})} \sigma_{I \setminus
  \{i\},I}.
\end{equation}
Then $\phi_i$, $i \in S$ anticommute with each other and with the
differential $d$, so that the homology $H_\idot(\kk{V})$
becomes a module over the exterior algebra $\Lambda^\hdot\langle
\phi_1,\dots,\phi_n\rangle$ generated by $\phi_1,\dots,\phi_n$.

\begin{defn}\label{polytight}
A $G_S$-module $V$ is called {\em tight} if the homology
$H_\idot(\kk{V})$ is the free $\Lambda\langle \phi_1,\dots,\phi_n
\rangle$-module generated by $H_0(\kk{V})$.
\end{defn}

It is easy to see that when $n=|S|=1$, this reduces to
Definition~\ref{tight}, so that the terminology and notation is
consistent with Subsection~\ref{tight.sub}.

We will now give a different construction of the complex $\kk{V}$
which is more canonical (and in particular, manifestly independent
of the choice of an order on $S$). Namely, let $\C'$ be the category
of $G_S$-modules in $\C$. For any $i \in S$, let $R_i = k[\Z/p\Z]$
be the regular representation of $\Z/p\Z$ over the field $k$
considered as a $G_S$-module by means of the natural projection $G_S
\to G_{\{i\}}=\Z/p\Z$. Define a functor $F:\C' \to \C'$ by
$$
F(V) = \bigoplus_{i \in S}V \otimes R_i.
$$
This is obviously an exact functor; the trace projections $R_i \to
k$ induce a map $F \to \Id$. Therefore we can set up a bar
resolution as described in Subsection~\ref{prelim} and obtain a
functorial resolution $B_\idot(F)(V)$. We form a canonical complex
$\cc{V}$ by setting
\begin{equation}\label{cc.def}
\cc{V}_\idot = H^0(G_S,B_\idot(F)(V)).
\end{equation}
Explicitly, we have
$$
B_l(F)(V) = V \otimes \left(\bigoplus R_i\right)^{\otimes l} = V
\otimes \bigoplus_{i_1,\dots,i_l \in S} R_{i_1} \otimes \dots
\otimes R_{i_l},
$$
and this space is naturally graded by the number of different
multiples on the right-hand side: we set
\begin{equation}\label{gr.B}
B_{l-m,m}(V) = V \otimes \bigoplus_{|\{i_1,\dots,i_l\}|=m} R_{i_1}
    \otimes \dots \otimes R_{i_l}.
\end{equation}
It is easy to see that this is a bicomplex -- the differential $d$
splits into a sum of a component $d^{1,0}$ of bidegree $(1,0)$ and a
component $d^{0,1}$ of bidegree $(0,1)$. This grading and the
bicomplex structure descend to the complex $\cc{V}$. The grading can
be further refined by specifying precisely the subset $\{
i_1,\dots,i_l\} \subset S$ -- for any $I \subset S$, we set
$$
B_{l-|I|,I}(V) = V \otimes \bigoplus_{\{i_1,\dots,i_l\}=I} R_{i_1}
\otimes \dots \otimes R_{i_l}.
$$

\begin{lemma}\label{bar=>dg}
For every $I \subset S$, the complex $B_{\idot,I}(V)$ is a free
resolution of $V$ considered as a $G_I$-module. Moreover, we have
$$
H_0(\cc{V},d^{0,1}) \cong \kk{V}
$$
-- in other words, $\kk{V}$ is the $0$-th homology complex of the
bicomplex $\cc{V}$ with respect to the differential $d^{0,1}$.
\end{lemma}

\proof{} Every product
$$
R=R_{i_1} \otimes \dots \otimes R_{i_l}
$$
is obviously a free $G_I$-module, where $I=\{i_1,\dots,i_l\}$, while
$G_{\overline{I}}$ acts trivially on $R$. Therefore the complex
$B_{\idot,I}(V)$ indeed consists of free $G_I$-modules, and
moreover,
$$
H^0(G_S,B_{\idot,I}(V)) \cong
H^0(G_I,B_{\idot,I}(H^0(G_{\overline{I}},V))) \cong
H_0(G_I,B_{\idot,I}(H^0(G_{\overline{I}},V))).
$$
Thus the second claim follows from the first, and to prove the
first, it suffices to prove that the complex $B_{\idot,I}(V)$ is
indeed a resolution of $V$. It obviously suffices to consider $\C =
k\mmod$, $V=k$ with the trivial $G_S$-action. By induction on
$n=|S|$, it suffices to consider $I=S$, and we may assume the claim
proved for all proper subsets $I \subset S$. Then instead of proving
that $B_{\idot,n} \to k$ is a quasiisomorphism, we may prove that
$B_\idot \to k_\idot$ is a quasiisomorphism, where $k_\idot$ is the
complex given by
$$
k_m = \bigoplus_{|I|=m}k, \qquad d = \sum_{i \in \overline{I}}
(-1)^{l(i,I)} \id.
$$
This latter complex is easily seen to be acyclic, while the
bar-complex $B_\idot$ is acyclic by construction. This finishes the
proof.
\endproof

We see that the order on $S$ that we used to construct the complex
$\kk{V}$ is an artefact of the particular explicit construction: the
bar resolution does not depend on the order and gives the same
thing. To see the maps $\phi_i$ in this language, we note that the
subgroup of $G_S$-invariants in $R_i$ is the trivial $G_S$-module
$\Z$; this gives a canonical embedding $k \to R_i$ and a functorial
map $\phi_i:\Id \to F$. As in \eqref{boott}, this gives rise to
functorial maps
$$
\phi_i:B_\idot(V)[1] \to B_\idot, \qquad \phi_i:\cc{V} \to \cc{V}
$$
and on the level of homology $H_0(\cc{V},d^{0,1})$, this is the
same map $\phi_i$ as in \eqref{phi.i}. We can also collect all the
maps $\phi_i$ together into a single map
\begin{equation}\label{phi.grp}
\phi:k[S] \otimes B_\idot(V)[1] \to B_\idot(V).
\end{equation}
One final result we will need is the following. Consider the
multiplicative group $(\Z/p\Z)^*$, and let it act on $G_S=\Z/p\Z[S]$
by dilations. Assume that a $G_S$-module $V$ is tight, and that the
$G_S$-action on $V$ is extended to an action of the semidirect
product $G_S \rtimes (\Z/p\Z)^*$. Moreover, assume that the induced
$(\Z/p\Z)^*$-action on the homology $H_0(\kk{V})$ is trivial.

\begin{lemma}\label{p-2.loc}
In the assumptions above, the canonical map
$$
\tau_{[0,l]}\cc{V}^{(\Z/p\Z)^*} \to \tau_{[0,l]}\kk{V}
$$
is a quasiisomorphism whenever $l \leq 2(p-2)$.
\end{lemma}

\proof{} The canonical maps $\phi_i:H_\idot(\kk{V}) \to
H_{\idot+1}(\kk{V})$ are obviously $(\Z/p\Z)^*$-equivariant;
therefore in the assumptions above, $(\Z/p\Z)^*$ acts trivially on
the homology $H_l(\kk{V})$ for any $l$. By Lemma~\ref{bar=>dg}, it
suffices to prove that for every $I \subset S$, the natural map
$$
H_\idot(G_I,H^0(G_{\overline{I}},V))^{(\Z/p\Z)^*} \to V_I
$$
is a quasiisomorphism in degrees $\leq 2(p-2)$ -- in other words, we
have to check that $H_l(G_I,H^0(G_{\overline{I}},V))$ has no
$(\Z/p\Z)^*$-invariant elements when $0 < l \leq 2(p-2)$. But
indeed, by the same assumption of tightness, we have
$$
H_l(G_I,H^0(G_{\overline{I}},V)) \cong V_I \otimes H_l(G_I,k),
$$
and $(\Z/p\Z)^*$ acts through the second factor on the right-hand
side, so that it suffices to check that
$H^0((\Z/p\Z)^*,H_l(G_I,k))=0$ for $0 < l \leq 2(p-1)$. This is
obvious: the homology space $H_\idot(G_I,k)$ is the free
module of rank $1$ over the free graded-commutative algebra -- or
rather, coalgebra -- $k[v_1,\dots,v_m,\phi_1,\dots,\phi_m]$
generated by $m=|I|$ elements $v_1,\dots,v_m$ of degree $-2$ and $m$
elements $\phi_1,\dots,\phi_m$ of degree $-1$, and all the
generators are $(\Z/p\Z)^*$-eigenvectors with weight $1$.
\endproof

\begin{remark}
It is rather unfortunate that we have to define the complex $\cc{V}$
as an explicit complex, by fixing an explicit resolution
$B_\idot(V)$. It would be much nicer to be able to define it as a
derived functor of some type. Then $|S|=1$, it is possible, and we
indeed do this in Lemma~\ref{tate}: the complex $\cc{V}$ can be
equivalently defined as the negative part of the canonical
filtration in the Tate homology complex $H_{\#}(G_S,V)$; if $V$ is
tight, the homology of the complex $\cc{V}$ is isomorphic to
$$
H_{\#}(G_S,V)/H^\hdot(G_S,V) = H^\hdot(G_S,V) \otimes_{k[[u]]}
k[u^{-1}],
$$
where $u \in H^2(G_S,\Z)$ is the periodicity generator, and none of
this depends on any choices of a free resolution. When $n=|S| > 1$,
we still have
$$
H_\idot(\cc{V}) \cong H^\hdot(G_S,V) \otimes_{k[[u_1,\dots,u_n]]}
k[u_1^{-1},\dots,u_n^{-1}],
$$
but this is not a good description from the general homological
point of view (in particular, it makes no sense in the derived
category and it is not sufficiently functorial in $S$). A derived
category interpretation is unknown, and none of the obvious
candidates give the correct answer. A similar situation occurs, for
instance, in the representation theory of affine Kac-Moody algebras
(see e.g. \cite{FF}): one can construct a certain type of homology
by an explicit resolution, and while it is intuitively clear that
the choice of a resolution should be irrelevant, there is no general
categorical framework which makes it precise.
\end{remark}

\subsection{Polycyclic categories.}
For any set $X$, define a natural contravariant functor
$X^{\#}:\Lambda^{opp} \to \Sets$ by
\begin{equation}\label{hash}
X^{\#}([m]) = \Maps(\Lambda([1],[m]),X) = X^m.
\end{equation}
If $X = G$ is a group, then $G^{\#}$ can be treated as a functor
from $\Lambda^{opp}$ to the category of groups.

\begin{defn}
The {\em wreath product} $G \wrth \Lambda$ is the fibered category
over $\Lambda$ which is the total space of the functor $G^{\#}$.
The {\em polycyclic categories} $\B_p$, $p \geq 1$, and $\B_\infty$
are the wreath products $\B_p = (\Z/p\Z) \wrth \Lambda$, $\B_\infty
= \Z \wrth \Lambda$.
\end{defn}

Explicitly, objects of $G \wrth \Lambda$ are $[n]$, $n \geq 1$, and
for any $[n],[m]$ the set of morphisms from $[m]$ to $[n]$ is given
by
$$
\begin{aligned}
(G \wrth \Lambda)([m],[n]) &= G^m \times \Lambda([m],[n]) \\
& = \left\{ \langle f',f \rangle \mid f' \in \Hom^\sigma([m],G), f
\in \Hom_{\Cycl}([m],[n])/\sigma \right\},
\end{aligned}
$$
with composition defined by $\langle g',g \rangle \circ \langle f',
f \rangle = \langle (g' \circ f)f', g \circ f \rangle$. Here we
recall that $\Cycl$ is the category of linearly ordered sets
equipped with an order-preserving automorphism $\sigma$, for any $l
\geq 1$, $[l] \in \Cycl$ is the linearly ordered set $\Z$ with
$\sigma:x \mapsto x + l$, and we interpret $G^m$ as the set of
$\sigma$-invariant maps from $[m] = \Z$ to $G$.

In the particular case $G = \Z$, $G \wrth \Lambda = \B_\infty$, we
define for any $[m]$, $[n]$ a subset $\overline{\B}_\infty([m],[n])
\subset \B_\infty([m],[n])$ by
\begin{itemize}
\item $\langle f',f \rangle \in \overline{\B}_\infty([m],[n])$ if
  and only if for any $l$, $0 \leq l < m$, we have 
\begin{equation}\label{remndr}
0 \leq f(l) +  nf'(l) < n.
\end{equation}
\end{itemize}
It is easy to see that \eqref{remndr} is preserved by the
composition law, so that it defines a subcategory
$\overline{B}_\infty \subset \B_\infty$. Moreover, for any $f \in
\Hom_{\Cycl}([m],[n])$ there is exactly one $f'$ such that $\langle
f',f \rangle$ satisfy \eqref{remndr}. We conclude that
$\overline{\B}_\infty \cong \Lambda_\infty$, so that we have a
canonical embedding $\lambda_\infty:\Lambda_\infty \to
\B_\infty$. Reducing this modulo $p$, we obtain a canonical
embedding $\lambda_p:\Lambda_p \to \B_p$.

For any small category $\Sigma$ equipped with a functor $\Sigma \to
\Lambda$, we denote $G \wrth \Sigma = \Sigma \times_{\Lambda} G
\wrth \Lambda$. In particular, we have the simplicial category
$\Delta^o$ equipped with the functor $j:\Delta^o \to \Lambda$, and
we can form the wreath product $G \wrth \Delta^o$. Denote
$\Delta^o_p = (\Z/p\Z) \wrth \Delta^o$. We have a natural embedding
$\Delta^o_p \to \B_p$, and a Cartesian diagram
\begin{equation}\label{D.p.0}
\begin{CD}
\Delta^o \times \ppt_p @>{j_p}>> \Lambda_p\\
@VVV @VV{\lambda_p}V\\
\Delta^o_p @>>> \B_p.
\end{CD}
\end{equation}
But by definition, for any $[m] \in \Delta^o$, the set
$\Lambda([1],j([m]))$ has a distinguished element; therefore for any
set $X$ the pullback $j^*X^{\#}$ admits a canonical projection
$j^*X^{\#} \to X_{\Delta^o}$ onto the constant functor
$X_{\Delta^o}:\Delta^o \to \Sets$. If $X = \Z/p\Z$, then this is
compatible with the group structure and induces a projection
$\Delta^o_p \to \Delta^o \times \ppt_p$. Denote
$$
\overline{\Delta^o}_p = \Delta^o_p \times_{\Delta^o \times \ppt_p}
\Delta^o.
$$
Then the Cartesian diagram \eqref{D.p.0} extends to a diagram
\begin{equation}\label{D.p}
\begin{CD}
\Delta^o @>>> \Delta^o \times \ppt_p @>{j_p}>> \Lambda_p\\
@VVV @VVV @VV{\lambda_p}V\\
\overline{\Delta^o}_p @>>> \Delta^o_p @>>> \B_p
\end{CD}
\end{equation}
with Cartesian squares.  By abuse of notation, we will denote the
composition $\overline{\Delta^o}_p \to \Delta^o_p \to \B_p$ of the
embeddings in the bottom row by the same letter $j_p$.

Finally, note that the multiplicative group $(\Z/p\Z)^*$ acts on
$\Z/p\Z$ and consequently on $\B_p$, so that we can form the
semidirect product $\wt{\B}_p=\B_p \rtimes (\Z/p\Z)^*$: it has the
same objects $[n]$, $n \geq 1$, we set
$$
\wt{\B}_p([m],[n]) = \left\{ \langle f,f' \rangle | f \in
\B_p([m],[n]), f' \in (\Z/p\Z)^*\right\},
$$
and the composition law is $\langle g,g' \rangle \circ \langle f,f'
\rangle = \langle f'(g)f,g'f' \rangle$. We have a natural embedding
$\B_p \subset \wt{\B}_p$, and the embedding $\lambda_p:\Lambda_p \to
\B_p$ extends to an embedding $\wt{\lambda}_p:\Lambda_p \to
\wt{\B}_p$. We will call $\Fun(\wt{\B}_p,\C)$ {\em extended
polycyclic category}, and we will call its objects {\em extended
polycyclic}.

\subsection{Tight polycyclic objects.}

For any $E \in \Fun(\B_p,\C)$, we can form the bar resolution
$B_\idot^!(\lambda_p^*E) = B_\idot(j_{p!}j_p^*)(\lambda_p^*E) \in
\Fun(\Lambda_p,\C)$ as in Subsection~\ref{bar.sub}. Applying the
base change to the diagram \eqref{D.p}, we see that the functor
$j_{p!}j_p^*$ extends to a functor $j_{p!}j_p^*:\Fun(\B_p,\C) \to
\Fun(\B_p,\C)$, so that the resolution $B_\idot^!(\Lambda_p^*E)$
comes from a polycyclic resolution $B_\idot^!(E) \in
\Fun(\B_p,\C)$. For any $[n] \in \B_p$, the complex
$B_\idot^!(E)([n])$ is the functorial bar resolution
$B_\idot(E)(E([n]))$ of the polycyclic module $E([n])$ that we have
considered in Subsection~\ref{poly.grp}. If we apply the direct
image functor $\chi_*:\Fun(\B_p,\C) \to \Fun(\Lambda,\C)$, then we
have
$$
\chi_*(B_\idot^!(E))([n]) \cong \cc{E([n])}_\idot,
$$
where $\cc{E([n])}_\idot$ is as in \eqref{cc.def}. The additional
grading \eqref{gr.B} on $B_\idot^!(E)([n])$ is not compatible with
maps $[n] \to [n']$ in the category $\Lambda$. However, one checks
easily that the associated filtration is preserved by the maps: if
we set
$$
F_qB_l^!(E)([n]) = \bigoplus_{m \geq q}B_{m,l-m}(E([n]))
$$
for every $[n]$, then these fit together into a subobject
$F_qB_l^!(E) \subset B_l^!(E)$. Denote
$$
\kk{E} = \chi_*(B_\idot^!(E))/(\chi_*(F_1B_\idot^!(E)) + d
\chi_*(F_1B_\idot^!(E))),
$$
where $d$ is the differential in the complex $B_\idot^!(E)$. Then
$\kk{E}$ is a well-defined complex in $\Fun(\Lambda,\C)$, functorial
in $E$, and Lemma~\ref{bar=>dg} immediately shows that we have
$$
\kk{E}([n]) = \kk{E([n])}
$$
for any $[n] \in \Lambda$.

Recall that we have defined in Subsection~\ref{bar.sub} a canonical
map $\pi^*(j_!k) \to j_{p!}k$ of $p$-cyclic objects; the reader will
check easily that this map is compatible with the polycyclic
stucture on $j_{p!}k$ and gives a map $\phi:\chi^*(j_!k) \to
j_{p!}k$ of objects in $\Fun(\B_p,k\mmod)$. Composing this with
\eqref{boott}, we obtain a functorial map
$$
\phi:\chi^*j_!k \otimes B_\idot^!(E)[1] \to B_\idot^!(E)
$$
for any $E \in \Fun(\B_p,\C)$, and as in Subsection~\ref{bar.sub},
this allows to introduce the conjugate filtration $W_\idot$ on
$B_\idot^!(E)$ and on $\chi_*(B_\idot^!(E)$. We note that if we
equip $\kk{E}$ with the canonical filtration, then the natural map
$\chi_*(B_\idot^!(E)) \to \kk{E}$ is a filtered map.  By the
projection formula, the map $\phi$ induces a map
$$
\phi:j_!k \otimes \chi_*B_\idot^!(E)[1] \to \chi_*B_\idot^!(E).
$$
We note that for any $[n] \in \Lambda$, we have $j_!k([n]) = k[S]$,
where the finite set $S$ is the set of maps
$\Maps_{\Lambda}([1],[n])$, and $\phi$ evaluated at $[n]$ is just
the canonical map \eqref{phi.grp}. In particular, $\phi$ induces a
map
$$
\phi:j_!k \otimes \kk{E}[1] \to \kk{E}.
$$
Restricting this to $k(1) \subset j_!k$, we obtain a map
$$
v:\kk{E}(1)[1] \to \kk{E}.
$$
It turns out that a statement completely analogous to
Lemma~\ref{adag}, \eqref{tri.1} is true for the complex $\kk{E}$.
Namely, for any $E \in \Fun(\B_p,\C)$ and any $[n] \in \B_p$, the
object $E([n]) \in \C$ by definition carries a representation of the
polycyclic group $(\Z/p\Z)^n$ -- or, more precisely, of the group
$\Z/p\Z[S]$, where $S = \Maps([1],[n])$.

\begin{defn}\label{polytight.cat}
The object $E \in \Fun(\B_p,\C)$ is called {\em tight} if
\begin{enumerate}
\item for any object $[n] \in \B_p$, the $(\Z/p\Z)^n$-module
$E([n])$ is tight in the sense of Definition~\ref{polytight}, and,
\item the object $\lambda_p^*E
\in \Fun(\Lambda_p,\C)$ is tight in the sense of
Definition~\ref{tight.cat}, and the natural map
$$
\lambda_p^*:\ii{E} \to \ii{\lambda_p^*E}
$$
is an isomorphism (where we denote $\ii{E} = H_0(\kk{E})$).
\end{enumerate}
\end{defn}

\begin{remark}
The condition~\thetag{ii} of Definition~\ref{polytight.cat} is
probably redundant, but were not able to prove it, and in practical
applications, it is very easy to check the condition by hand.
\end{remark}

\begin{lemma}\label{kk.E}
Let $E \in \Fun(\B_p,\C)$ be tight in the sense of
Definition~\ref{polytight.cat}. Then there exists a triangle
$$
\begin{CD}
\kk{E}(1)[1] @>{v}>> \kk{E} @>>> \kk{\lambda_p^*E} @>>>
\end{CD}
$$
which, if we equip all the terms with the canonical filtration,
becomes a filtered exact triangle after applying the cyclic homology
functor $HC_\idot$. 
\end{lemma}

\proof{} By definition, we have to prove that for any $l \geq 0$,
the cyclic homology functor $HC_\idot$ turns the sequence
\begin{equation}\label{sq}
\begin{CD}
0 @>>> \kkk{E}{l-1}(1) @>{v}>> \kkk{E}{l} @>>> \kkk{\lambda_p^*E}{l}
@>>> 0
\end{CD}
\end{equation}
into an exact triangle. By Corollary~\ref{per.corr}, we may equally
well use the Hochschild homology functor $HH_\idot$.  For $l=0$, the
sequence is exact already in $\Fun(\Lambda,\C)$ by
Definition~\ref{polytight.cat}~\thetag{ii}. By
Definition~\ref{polytight}, we have
$$
\kkk{E}{l} \cong \ii{E} \otimes \Lambda^lj_!k,
$$
and the map $v$ is induced by the exterior multiplication map
$m:j_!k \otimes \Lambda^{l-1}j_!k \to \Lambda^lj_1k$ restricted to
$k(1) \subset j_!k$ (here, sadly, $\Lambda$ with an upper index has
to mean exterior power in the sense of the pointwise tensor product
in $\Fun(\Lambda,k\mmod)$, and we hope that this does not cause any
confusion). By Definition~\ref{tight}, $\kkk{\lambda_p^*E}{l}$ is
isomorphic to $\ii{E}$ for $l=1$ and trivial for $l \geq
2$. Consider the multiplication map
\begin{equation}\label{m.eq}
m:\Lambda^{l-1}j_1k \otimes k(1) \to \Lambda^lj_!k.
\end{equation}
If $l=1$, then this is just the embedding map $k(1) \to j_!k$; it is
injective, its cokernel is by definition the constant functor $k$,
and we conclude that \eqref{sq} is exact for $l=1$, also already in
$\Fun(\Lambda,\C)$. For $l \geq 2$, we restrict terms to
$\Fun(\Delta^o,k\mmod)$ and use the Dold-Thom equivalence $\DT$
\cite{DT} between simplicial objects in an abelian category and
non-positively graded complexes in the same category. Under this
equivalence, $\DT(j^*k(1)) \cong k[1]$, $\DT(j_!k) \cong k[1] \oplus
k$, and moreover, $\DT(\Lambda^lj_!k) \cong k[l] \oplus
k[l-1]$. Therefore \eqref{m.eq} would be an isomorphism, were the
functor $\DT$ compatible with the tensor product. Since it is only
compatible with the tensor product up to a quasiisomorphism, we
conclude that $\DT(m)$ is a quasiisomorphism. Therefore for any $E'
\in \Fun(\Delta^o,\C)$, in particular for $E' = j^*\ii{E}$, the map
$$
\DT(m):\DT(E' \otimes \Lambda^{l-1}j_1k \otimes k(1)) \to \DT(E'
\otimes \Lambda^lj_!k)
$$
is a quasiisomorphism. Since for any $E' \in \Fun(\Delta^o,\C)$,
$H_\idot(\Delta^o,E')$ is the homology of the complex $\DT(E')$, we
conclude that $v$ in \eqref{sq} becomes a quasiisomorphism after
applying $H_\idot(\Delta^o,-) = HH_\idot(-)$, so that \eqref{sq}
indeed becomes exact.
\endproof

\subsection{Comparison maps.}
Below in Section~\ref{alg}, in our applications to cyclic homology
for associative algebras, the computation at some point passes
through a tight polycyclic object $E$, and it turns out that the
associated complex $\kk{E}$ is rather easy to control. However, it
is the complex $\pi_{\flat}\lambda_p^*E$ that is related to the
cyclic homology. Both $HC_\idot(\kk{E})$ and
$HC_\idot(\pi_{\flat}\lambda_p^*)$ have very similar structure -- in
particular, both are equipped with filtrations whose associated
graded quotients are the same (and isomorphic to
$HH_\idot(\ii{E})[v]$). It would be very nice to know that the two
complexes themselves are quasiisomorphic, possibly under additional
natural assumptions on $E$. Unfortunately, we were not able to prove
it, and there are reasons to believe that it is not true: the
extension data between associated graded pieces of the corresponding
filtrations {\em are} different. In particular, there is no natural
map from one complex to the other. We had to settle for a weaker
comparison result: there is a third complex which maps both into
$\kk{E}$ and into $\pi_{\flat}\lambda_p^*E$, and after we take
cyclic homology, both maps become quasiisomorphisms {\em in degrees
$\geq -2(p-2)$}. Consequently, $\kk{E}$ and $\pi_{\flat}\lambda_p^*$
do have equal cyclic homology in low degrees, which turns out to be
enough in application to the Hodge-to-de Rham degeneration. An
instance of this phenomenon occurs already in the commutative
situation considered in \cite{DL} (where the characteristic of the
base field had to be greater than the dimension of the algebraic
variety in question).

The comparison result that we prove is as follows. Assume given a
polycyclic object $E \in \Fun(\B_p,\C)$ which is tight in the sense
of Definition~\ref{polytight.cat}. Moreover, assume that $E$ comes
from an {\em extended} polycyclic object $E \in
\Fun(\wt{\B}_p,\C)$. It is easy to check that the bar resolution
$B_\idot^!(E)$ together with its conjugate filtration $W_\idot$ is
compatible with the extended polycyclic structure, so that we have a
natural comparison map
\begin{equation}\label{comp.1}
\wt{\chi}_*(B_\idot^!(E)) \to \chi_*(B_\idot^!(E)) \to \kk{E},
\end{equation}
which becomes a map of filtered complexes in $\Fun(\Lambda,\C)$ if
we equip $\kk{E}$ with the canonical filtration. On the other hand,
we have a natural filtered map
\begin{equation}\label{comp.2}
\wt{\chi}_*(B_\idot^!(E)) \to \pi_*(B_\idot^!(\lambda_p^*E)).
\end{equation}
Let us say that a map $f:E_\idot \to E_\idot'$ of (filtered)
complexes in $\Fun(\Lambda,\C)$ is a (filtered) quasiisomorphism
{\em up to homology} in degrees $\geq m$ if the correponding map
$f:HC_\idot(E_\idot) \to HC_\idot(E_\idot')$ is a (filtered)
quasiisomorphism in degrees $\geq m$.

\begin{prop}\label{trunc}
In the assumptions above, the comparison maps \eqref{comp.1} and
\eqref{comp.2} are filtered quasiisomorphisms up to homology in
degrees $\geq \!-2(p-2)$.
\end{prop}

\proof{} The map \eqref{comp.1} is a quasiisomorphism in degrees
$\geq -2(p-2)$ by Lemma~\ref{p-2.loc}, and since it is compatible
with the periodicity map, it is a filtered quasiisomorphism up to
homology in the said degrees by Lemma~\ref{kk.E}. The map
\eqref{comp.2} is also compatible with the periodicity map, and by
the definition of the conjugate filtration, it suffices to prove
that the induced map
$$
\gr^W_1\wt{\chi}_*(B_\idot^!(E)) \to
\gr^W_1\pi_*(B_\idot^!(\lambda^*_pE))
$$
is a quasiisomorphism up to homology in degrees $\geq -2(p-2)$. By
definition, both sides are complexes concentrated in degrees $\leq
-1$. In degree $-1$, the homology of the left-hand side isomorphic
to $j_!j^*\ii{E}$, while the homology of the right-hand side is
$j_!j^*\ii{\lambda_p^*E}$, and the map is an isomorphism by
Definition~\ref{polytight.cat}~\thetag{ii}; thus it suffices to
prove that $\tau_{\leq -2}\gr^W_1\wt{\chi}_*(B_\idot^!(E))$ and
$\tau_{\leq -2}\gr^W_1\pi_*(B_\idot^!(\lambda^*_pE))$ are trivial up
to homology in degrees $\geq -2(p-2)$. For the first complex, we can
apply the map \eqref{comp.1} and deduce the statement. For the
second complex, we note that for any tight $p$-cyclic $E' \in
\Fun(\Lambda_p,\C)$, by Lemma~\ref{bar.perio} and
Corollary~\ref{bar.perio.corr} $\tau_{\leq
-2}\gr^W_1\pi_*B_\idot^!(E')$ is trivial up to homology in all
degrees.
\endproof

\section{Associative algebras.}\label{alg}

We now fix a perfect field $k$ of odd positive characteristic $p$,
and we assume given a $k$-linear Grothendieck abelian category
$\C$. Moreover, we assume that $\C$ is equipped with a symmetric
tensor product with unit object $k \in \C$.

\subsection{Definitions.}\label{alg.def} 
Let $A$ be an associative unital algebra object in the tensor
category $\C$. We define a canonical cyclic object $A_{\#} \in
\Fun(\Lambda,\C)$ in the following way. For any $m \geq 1$, we set
$$
A_{\#}([m])=A^{\otimes m} = \bigotimes_{i \in \overline{[m]}}A,
$$
where we number the factors in the tensor product by elements of the
set $\overline{[m]} = \Lambda([1],[m])$. For any map $f \in
\Lambda([m],[n])$, we set
$$
f_{\#} = \bigotimes_{i \in \overline{[n]}}f_i,
$$
where
$$
f_i:\bigotimes_{j \in \overline{[f]}^{-1}(i)}A \cong A^{\otimes
  |\overline{f}^{-1}(i)|} \to A
$$
is the multiplication map given by the algebra structure on $A$, and
$\overline{f}:\overline{[m]} \to \overline{[m]}$ is the natural map
induced by $f \in \Lambda([m],[n])$ (if $\overline{f}^{-1}(i)$ is
empty, we let $A^{\otimes 0} = k$ be the unit object in $\C$, and we
take $f^i:k \to A$ to be the unit embedding).

This is a well-defined cyclic object in $\C$, so that we have the
Hochschild homology complex $HH_\idot(A_{\#}) \in D^-(\C)$ and the
cyclic homology complex $HC_\idot(A_{\#}) \in D^-(\C)$.

To proceed further, we specialize to the following situation. Assume
fixed a small category $Z$ equipped with a Grothendieck topology
$J$. From now on, and until Section~\ref{main}, let $\C$ be the
category of sheaves of $k$-vector spaces on $\langle Z,J
\rangle$. Moreover, assume given and fixed a cohomological functor
$H^\hdot(-)$ from $\C$ to the category of $k$-vector spaces (for
example, this may be the cohomology of the site $\langle Z,J
\rangle$, or the cohomology with some fixed supports).

\begin{defn}
The {\em Hochschild and cyclic homology} of the algebra $A \in \C$
is given by
$$
HH_\idot(A) = H^\hdot(\langle Z,J \rangle, HH_\idot(A_{\#})) \qquad
HC_\idot(A) = H^\hdot(\langle Z,J \rangle, HC_\idot(A_{\#})),
$$
where $HH_\idot(A_{\#})$ and $HC_\idot(A_{\#})$ are as in
\eqref{hh.hc}.
\end{defn}

The Hodge and the conjugate filtrations on $HC_\idot(A_{\#})$ induce
filtrations on $HC_\idot(A)$. We will be interested in the conjugate
filtration $W_\idot$. By abuse of notation, for any $l \geq 0$ we
will denote
$$
\begin{aligned}
W_lHC_\idot(A) &= H^\hdot(\langle Z,J \rangle,
W_l\pi_!i_p^*A_{\#}),\\
\gr^W_lHC_\idot(A) &= H^\hdot(\langle Z,J \rangle, \gr^W_l\pi_!i_p^*A_{\#}),
\end{aligned}
$$
although the map $W_lHC_\idot(A) \to HC_\idot(A)$ does not have to
be injective, so that $W_\idot HC_\idot(A)$ is only a filtration in
the generalized sense.

For any $V \in \C$, denote by $V^\tw = \Fr^*(V)$ the sheaf
$V$ with $k$-vector space structure twisted by the Frobenius map
$\Fr:k \to k$.

\begin{lemma}\label{otimesp}
For any $V \in \C$, the tensor power $V^{\otimes p}$ equipped with
the natural action of $\Z/p\Z$ by transpositions is tight in the
sense of Definition~\ref{tight}, and there is a canonical
isomorphism $\ii{V^{\otimes p}} \cong V^\tw$.
\end{lemma}

\proof{} Since taking the associated sheaf is an exact functor, it
suffices to prove the statement for the trivial topology $J$ on $Z$,
so that $\C$ is the category of presheaves of $k$-vector spaces on
$Z$. This in turn reduces to proving the statement for the vector
spaces $V(X)$ for all objects $X \in Z$. Thus we may assume that $Z$
is trivial and $V = k[S]$ is just a $k$-vector space with some basis
$S$. Then $V^{\otimes p} = k[S^p] = k[S] \oplus V'$, where $S
\subset S^p$ is the diagonal, and $V' = k[S^p \setminus S]$ is the
natural complement to $V = k[S] \subset V^{\otimes p}$. The group
$\Z/p\Z$ acts freely on $S^p \setminus S$; therefore $V'$ is a free
$\Z/p\Z$-module, it is trivially tight, and we have $\ii{V'} =
0$. The action on $V \subset V^{\otimes p}$ is trivial, so that $V$
also tight, and $\ii{V}\cong V$. Therefore $V^{\otimes p} = V \oplus
V'$ is tight, and $\ii{V^{\otimes p}} \cong V$. Finally, to
construct a functorial isomorphism $V^\tw \cong \ii{V^{\otimes p}}$,
we note that the embedding $V \to V^{\otimes p}$ defined by the
basis $S$ sends $V$ into the subspace of $\Z/p\Z$-invariant vectors,
and moreover, coincides modulo $T(V^{\otimes p})$ with
$$
v \mapsto (v \otimes v \otimes \dots \otimes v) \mod T(V^{\otimes
  p}) \in V^{\otimes p}/T(V^{\otimes p}).
$$
This is a well-defined map, it is explicitly indepedent of any
bases, it is {\em a posteriori} additive, and the $k$-vector space
structure on both sides obviously differs by the Frobenius map
$\Fr:k \to k$, $\lambda \mapsto \lambda^p$.
\endproof

\begin{corr}\label{cor.tgt}
For any associative unital algebra $A$ in the tensor category $\C$,
the object $i_p^*A_{\#} \in \Fun(\Lambda_p,\C)$ is tight, and we
have $\ii{i_p^*A_{\#}} \cong A_{\#}^\tw$, so that
\begin{equation}\label{cartier}
\gr^W_l\Ll^\hdot\pi_!i_p^*A_{\#} \cong j^\dg A_{\#}^\tw[2l]
\end{equation}
for any $l \geq 1$, and $\gr^W_lHC_\idot(A) \cong HH_\idot(A^\tw)[2l]$.
\end{corr}

\proof{} By definition, for any $m \geq 1$ we have 
$$
i_p^*A_{\#}([m]) \cong A_{\#}([pm]) \cong
\left(A_{\#}([m])\right)^{\otimes p},
$$
this is tight by Lemma~\ref{otimesp}, and $\ii{i_p^*A_{\#}([m])}
\cong A_{\#}([m])^\tw$. We leave it to the reader to check that the
isomorphism is compatible with the action of maps $f \in
\Lambda([m],[n])$, $n,m \geq 1$. The last claim is \eqref{tri.2} of
Lemma~\ref{adag}.
\endproof

\begin{remark}
The isomorphism \eqref{cartier} and the induced isomorphism
$$
\gr^W_lHC_\idot(A) \cong HH_\idot(A^\tw)
$$
is probabely the closest analog of the usual Cartier isomorphism in
our non-commutative theory.
\end{remark}

\subsection{Polycyclic structure.} We now note that if the algebra $A$ is
equipped with an action of a group $G$, then the object $A_\# \in
\Fun(\Lambda,\C)$ has a natural structure of a functor from $G \wrth
\Lambda$ to $\C$. In particular, if $G = \Z/p\Z$, we have a
polycyclic object $A_\# \in \Fun(\B_p,\C)$. On the other hand, the
length-$2$ complex $\kk{A}$ has a natural structure of a DG algebra
in $\C$.

\begin{lemma}\label{dag}
In the assumptions above, we have
$$
\kk{A_\#} \cong \kk{A}_\#.
$$
\end{lemma}

\proof{} Evaluating on $[n] \in \B_p$, we have $A_\#([n]) =
A^{\otimes n}$, and the action of the group $(\Z/p\Z)^n$ on this
object is induced by the $(\Z/p\Z)$-action on each factor in
$A^{\otimes n} = A \otimes \dots \otimes A$. Therefore for any $I
\subset S$, $S = \Lambda([1],[n])$, we have
$$
A_\#([n])_I \cong \left(\bigotimes_{i \in I} H_0(\Z/p\Z,A)\right)
\otimes \left(\bigotimes_{i \not\in I} H^0(\Z/p\Z,A)\right).
$$
This identification is by definition compatible with the
differentials and gives an isomorphism $\kk{A_\#([n])} \cong
\kk{A}^{\otimes n}$. It remains to check that these isomorphisms are
compatible with the maps $[m] \to [n]$; this is easy and left to the
reader.
\endproof

In particular, for any associative unital algebra $A$ in $\C$, the
$p$-th power $A^{\otimes p}$ has a natural algebra structure, and
this algebra $A^{\otimes p}$ is acted upon by the cyclic group
$\Z/p\Z$ -- and in fact, by the whole symmetric group $S_p$ on $p$
letters. Therefore we can form a polycyclic object $A^{\otimes p}_\#
\in \Fun(\B_p,\C)$, and this object is in fact extended polycyclic.

\begin{lemma}\label{dag.kk}
The polycyclic object $A^{\otimes p}_\#$ is tight in the sense of
Definition~\ref{polytight.cat}; we have $\lambda_p^*A^{\otimes p}_\#
\cong i_p^*A_\#$ and $\ii{A^{\otimes p}_\#} \cong \ii{i_p^*A_\#}
\cong A_\#^\tw$.
\end{lemma}

\proof{} The isomorphism $\lambda_p^*A^{\otimes p}_\# \cong
i_p^*A_\#$ immediately follows from the definitions, and the
condition \thetag{i} of Definition~\ref{polytight.cat} is an
immediate corollary of Lemma~\ref{dag}. Moreover, by the same lemma
we have $\ii{A^{\otimes p}_\#([n])} \cong \ii{A^{\otimes
p}}^{\otimes n} \cong \ii{A^{\otimes p}}_\#([n])$; together with
Corollary~\ref{cor.tgt}, this gives an isomorphism $\ii{A^{\otimes
p}_\#} \cong \ii{i_p^*A_\#} \cong A_\#^\tw$, which yields
Definition~\ref{polytight.cat}~\thetag{ii}.
\endproof

\begin{defn}\label{dag.defn}
For any $V \in \C$, denote $V^\dg = \kk{V^{\otimes p}}$.
\end{defn}

\begin{lemma}\label{dag.dege}
Assume that for an associative algebra $A$ in $\C$, the DG algebra
$A^\dg$ is formal -- in other words, $A^\dg$ is quasiisomorphic to
its homology, the trivial square-zero extension of $A^\tw$ by
$A^\tw[1]$. Then the spectral sequence for the conjugate filtration
computing $W_1HC_\idot(A)$ degenerates up to the term $E^{p-2}$.
\end{lemma}

\proof{} Since the conjugate filtration is periodic, it suffices to
analyse the top terms of the spectral sequence -- that is, we have
to show that
$$
W_{[1,(p-2)]}HC_\idot(A) \cong \bigoplus_{1 \leq l \leq p-2}
\gr^W_lHC_\idot(A).
$$
In other words, we have to analyse the cyclic homology of
$W_{[1,(p-2)]}\pi_{\flat}i_p^*A_\#$. By Proposition~\ref{trunc}, we
may replace it with $\tau_{[-1,-(p-2)]}\kk{A^{\otimes p}_\#}$, and
by Lemma~\ref{dag.kk}, this is isomorphic to
$\tau_{[-1,-(p-2)]}A^\dg_\#$. In other words, we have to show that
the canonical filtration on the complex $A^\dg_\#$ splits in a
certain range of degrees. But by assumption, we have $A^\dg \cong A^\tw
\oplus A^\tw[1]$ -- which means that the complex $A^\dg_\#$ is
quasiisomorphic to the sum of its homology in all degrees.
\endproof

\subsection{Splittings.}\label{spl.sub}
We will now study the formality of the DG algebra $A^\dg$. To do
this, we need to refine (and explain) Lemma~\ref{otimesp}.

Let $V$ be a $k$-vector space. Denote by
\begin{equation}\label{rho.0}
\rho^0: V \to H^0(\Z/p\Z,V^{\otimes p})
\end{equation}
the map which sends $v$ to $v \otimes \dots \otimes v$. This map,
although not additive, is functorial with respect to $V$. Therefore
for any small category $Z$ and a presheaf $V \in \Fun(Z^o,k\mmod)$
of $k$-vector spaces on $Z$, we have a natural map $\rho:V \to
H^0(\Z/p\Z,V^{\otimes p})$. In particular, we may take $Z = \Delta$,
the category of linearly ordered finite sets. By the Dold-Thom
Theorem \cite{DT}, the category $\Fun(\Delta^o,k\mmod)$ is
equivalent to the category of non-positively graded complexes of
$k$-vector spaces. Take a $k$-vector space $V$, consider the complex
$V[1]$, denote the associated simplicial vector space by $V(1)$, and
consider the map
$$
\rho^1:V(1) \to H^0(\Z/p\Z,V(1)^{\otimes p}).
$$
This map extends to an {\em additive} map
\begin{equation}\label{rho.1}
\rho^1:k[V(1)] \to H^0(\Z/p\Z,V(1)^{\otimes p}).
\end{equation}
Apply now the Dold-Thom equivalence $\DT$. The left-hand side becomes
the standard bar complex $C_\idot(V,k)$ which computes the homology
of the vector space $V$ -- considered as an abelian group -- with
coefficients in the field $k$. On the other hand, one checks easily
that we have an isomorphism of $(\Z/p\Z)$-modules
$$
\DT_\idot(V(1)^{\otimes p}) \cong V^{\otimes p} \otimes
\DT_\idot(k(1)^{\otimes p}),
$$
and the complex $\DT(k(1)^{\otimes p})$ is quasiisomorphic to $k[p]$
and concentrated in degrees $\leq -1$, $\DT_1(k(1)^{\otimes p})$ is
the trivial $(\Z/p\Z)$-module $k$, and all the $\DT_l(k(1)^{\otimes
p})$, $l \geq 2$ are regular representations of
$(\Z/p\Z)$. Therefore we can take the standard periodic resolution
$I_\#(k)$ of the trivial $(\Z/p\Z)$-module $k$, set $I_\flat(k) =
\tau_{\leq 1}I_\#(k)$ as in Subsection~\ref{tate.sub}, and choose a
$(\Z/p\Z)$-module map $\psi:\DT(k(1)^{\otimes p}) \to I_\flat(k)[2]$
which is a quasiisomorphism in degrees $\geq -p$. Composing this
with $\DT(\rho)$, we obtain a canonical map
\begin{equation}\label{rho.V}
\begin{aligned}
C_\idot(V,k) &\to H^0(\Z/p\Z,DT_\idot(V(1)^{\otimes p})) \to\\
&\to H^0(\Z/p\Z,V^{\otimes p} \otimes I_\flat(k))[2] \cong
H^0(\Z/p\Z,I_\flat(V))[2].
\end{aligned}
\end{equation}
Note that that the only choice of the whole construction was the
choice of the map $\psi$, which does not depend on $V$ at
all. Therefore \eqref{rho.V} is completely functorial in $V$. In
particular, we may just as well let $V$ be a presheaf on some
category $Z$, and moreover, by applying the associated sheaf functor
we may even consider sheaves with respect to some topology $J$.

Truncating the right-hand side of \eqref{rho.V}, we obtain a
canonical map
\begin{equation}\label{cmpl}
C_\idot(V,k) \to V^\dg[1]
\end{equation}
which is again completely functorial in $V$. In degree $1$, this is
the canonical map $V \cong H_1(V,k) \to H_0(V^\dg)$ constructed in
Lemma~\ref{otimesp}. In degree $2$, we obtain a canonical cocycle
$\rho_V \in C^2(V,H_0(\Z/p\Z,V^{\otimes p}))$.

\begin{defn}
Denote by $\wt{V}$ the group obtained as an extension of $V$ by
$H_0(\Z/p\Z,V^{\otimes p})$ given by the cocycle $\rho$.
\end{defn}

One checks easily that the $2$-cocycle $\rho_V$ is in fact symmetric,
so that the group $\wt{V}$ is commutative. The group $\wt{V}$ is
{\em not} a $k$-vector space and not a group of characteristic $p$:
multiplication by $p$ is given by the map
$$
\wt{V} \to V \overset{\phi}{\longrightarrow} H_0(\Z/p\Z,V^{\otimes
  p}) \hookrightarrow \wt{V}.
$$
Assume from now on that the field $k$ is perfect; then one checks
easily that if we twist the $k$-module structure on $V$ by the
Frobenius map -- in other words, replace $V$ with $V^\tw$ -- then
the cocycle $\rho_V$ is compatible with multiplication by constants,
and $\wt{V}$ is in fact a module over the second Witt vectors ring
$W_2(k)$. Since the cocycle $\rho_V$ comes from a map of complexes
\eqref{cmpl}, it reduction modulo $V \subset H_0(\Z/p\Z,V^{\otimes
p})$ is the cocycle for the extension
$$
H^0(\Z/p\Z,V^{\otimes p}) \to V^\tw
$$
split as a map of sets by \eqref{rho.0}), so that we have $\wt{V}/p
\cong H^0(\Z/p\Z,V^{\otimes p})$.

\smallskip

To sum up: we have a natural three-step filtration 
$$
V^\tw \subset H_0(\Z/p\Z,V^{\otimes p}) \subset \wt{V}
$$
on the $W_2(k)$-module $\wt{V}$; the quotient $\wt{V}/V^\tw$ is
naturally identified with $H^0(\Z/p\Z,V^{\otimes p})$, and the
quotient $\wt{V}/H_0(Z/p\Z,V^{\otimes p})$ is naturally identified
with $V^\tw$. Multiplication by $p \in W_2(k)$ isomorphically sends
the top quotient $V^\tw$ of this filtration into the bottom
subobject $V^\tw \subset \wt{V}$.

\smallskip

The construction of the cocycle $\phi_V$ is also compatible with
tensor products in the following sense.

\begin{lemma}\label{prd}
Assume given two $k$-vector spaces $V$, $W$, and let $m$ be the
natural map
$$
H_0(\Z/p\Z,V^{\otimes p}) \otimes H^0(\Z/p\Z,W^{\otimes p}) \to
H_0(\Z/p\Z,(V \otimes W)^{\otimes p}).
$$
Then $m(\rho_V \otimes \rho^0_W) = \rho_{V \otimes W}$, where
$\rho^0_W$ is the map $\rho^0$ in \eqref{rho.0} for the vector space
$W$. Consequently, there exists a functorial map
$$
\wt{V} \otimes_{W_2(k)} \wt{W} \to \wt{V \otimes W}.
$$
\end{lemma}

\proof{} To obtain $\rho_V \otimes \rho^0_W$, one uses the same
procedure as for the cocycle $\rho_V$, but for the simplicial
abelian group $V(1) \otimes W$ instead of $V(1)$. The map $m$ is
induced by the natural map $V(1) \otimes W \to (V \otimes
W)(1)$. Since \eqref{rho.0} is functorial with respect to any maps
of simplicial groups, we get the desired compatibility.
\endproof

As before, all of the above obviously works not only for $k$-vector
spaces, but also for presheaves and sheaves of $k$-vector spaces on
a site $\langle Z,J \rangle$. Take a $k$-vector space or a sheaf
$V$. Denote by $\wt{V^\dg}$ the complex
$$
\begin{CD}
H_0(\Z/p\Z,V^{\otimes p}) @>{T}>> \wt{V}
\end{CD}
$$
placed in degrees $0$ and $-1$. This is a complex quasiisomorphic to
$V^\tw$, and we have a natural map $\wt{V^\dg} \to V^\dg$ which
induces an isomorphism on $H_0$. If $V=A$ is an associative unital
algebra, then by Lemma~\ref{prd} $\wt{A^\dg}$ is a DG algebra, and
the canonical map $\wt{A^\dg} \to A^\dg$ is a DG algebra map. Adding
the tautological embedding $A^\tw[1] \cong H_1(A^\dg)[1] \to A^\dg$,
we obtain a DG algebra quasiisomorphism
$$
\wt{A^\dg} \oplus A^\tw[1] \to A^\dg.
$$
Therefore the DG algebra $A^\dg$ is {\em always} quasiisomorphic to
the sum of tis homology, but in the wrong category: $\wt{A^\dg}$ is
a DG algebra over $W_2(k)$, not over $k$.

\begin{lemma}\label{formm}
Assume that there exists an associative flat $W_2(k)$-algebra
$\wt{A}$ in $\Shv(Z,J)$ such that $\wt{A}/p \cong A^\tw$. Then the DG
$k$-algebra $A^\dg$ is formal.
\end{lemma}

\proof{} We have two extensions of the algebra $A^\tw$: $\wt{A}$ and
$\wt{A^\dg}$.  Let $\overline{A^\dg}$ be their Baer difference --
that is, the middle cohomology of the complex
$$
\begin{CD}
A^\tw @>>> \wt{A} \oplus \wt{A^\dg} @>>> A^\tw,
\end{CD}
$$
where the left-hand side map is the sum, and the right-hand side map
is the difference of the natural maps. Then $\overline{A^\dg}$ is a
DG algebra over $W_2(k)$, and moreover, $p$ acts trivially on it, so
that $\overline{A^\dg}$ is in fact a DG algebra over $k$. On the
other hand, $\overline{A^\dg}$ is quasiisomorphic to $A^\tw$, and we
have a natural map $\overline{A^\dg} \to A^\dg$. Adding the
embedding $A^\tw[1] \to A^\dg$, we obtain the desired
quasiisomorphism $A^\tw \oplus A^\tw[1] \cong \overline{A^\dg}
\oplus A^\tw[1] \to A^\dg$.
\endproof

\begin{remark}\label{bok}
One can show without much difficulty that the converse to this
statement is also true: liftings of $V$ to a flat $W_2(k)$-module are
in functorial one-to-one correspondence with splittings of the
complex $V^\dg$ (understood in appropriate way), and this is
compatible with algebra structures, if they are present. We do not
go into this to save space.
\end{remark}

\begin{remark}\label{steen}
It is very interesting to repeat our construction of the map
$\rho^1$ in \eqref{rho.1} for different shifts -- one takes $V[l]$,
$l \geq 0$, instead of $V[1]$, and obtains a map $\rho^l$. If one
goes to the limit $l \mapsto \infty$, then the left-hand side of
\eqref{rho.1} becomes the {\em stable} homology of the group $V$ --
this is a complex which depends funtorially on $V$, and whose
homology is isomorphic to $V$ tensored with the dual to the Steenrod
algebra. The right-hand side stabilizes in a straightforward manner
and becomes quasiisomorphic to the truncated Tate homology $H_{\geq
-1}(\Z/p\Z,V^{\otimes p})$. The limit map $\rho^\infty$ is
essentially just the Steenrod $p$-th power map. The cocycle $\rho_V$
survives in the stable situation and becomes the Bokstein
homomorphism (this is the topological underlying reason for
Lemma~\ref{formm} and Remark~\ref{bok}). We note that the structure
of Steenrod algebra is well-known; in particular, after the Bokstein
homomorphism class in degree $1$, there is a non-trivial class in
degree $p$ (and in other higher degrees). We suspect that
$\rho^\infty$ sends these classes to some non-trivial classes in
$H_\idot(\Z/p\Z,V^{\otimes p})$, so that even if $V$ is lifted to a
$W_2(k)$-module, this Tate homology complex does not split as a
whole in a functorial way -- it only splits in degrees $\geq
-(p-1)$. Thus the restriction on degree in Lemma~\ref{dag.dege} and
Proposition~\ref{trunc} is unavoidable: this is how things really
are. We note that our trick of using the additional
$(\Z/p\Z)^*$-symmetry is lifted out the standard computation of
Steenrod powers found in any topology textbook.
\end{remark}

\section{Hodge to de Rham degeneration.}\label{main}

We are now ready to study the Hodge-to-de Rham spectral sequence and
its degeneration. Assume given a field $k$ and a small category $Z$
equipped with a Grothendieck topology $J$. Denote by $\C$ the
category of sheaves of $k$-vector spaces on $\langle Z,J
\rangle$. Assume also given a cohomological functor $H^\hdot(-)$
from $\C$ to the category of $k$-vector spaces.

We start with the positive characteristic case.

\begin{theorem}\label{char.p.main}
Assume given an associative unital algebra $A \in \C$ over a field
$k$ in the category $\C$ of sheaves of $k$-vector spaces on a site
$\langle Z,J \rangle$. Assume that $k$ is a perfect field of
positive odd characteristic $\cchar k=p > 2$. Moreover, assume that
\begin{enumerate}
\item the diagonal $A$-bimodule $A$ admits a finite flat resolution
  of length $n$,\label{reso}
\item for some integer $m > 0$ and every sheaf $E \in \C$, the
  cohomology groups $H^l(\langle Z,J \rangle,E)$ are trivial
  whenever $l > m$,
\item we have $p > n + m$, and\label{p.large}
\item for any integer $l$, the Hochschild homology group $HH_l(A)$
  is a finite-dimensional vector space over $k$.
\end{enumerate}
Finally, assume that there exists a flat algebra $\wt{A} \in
\Shv(Z,J)$ over the ring $W_2(k)$ of second Witt vectors of the
field $k$ such that $\wt{A}/p \cong A$.

Then the Hochschild-to-cyclic spectral sequence $HH_\idot(A)[v]
\Rightarrow HC_\idot(A)$ degenerates.
\end{theorem}

\begin{lemma}\label{conj.splits}
In the assumptions of Theorem~\ref{char.p.main}, the spectral
sequence associated to the conjugate filtration on $W_1HC_\idot(A)$
degenerates, so that we have an isomorphism $W_1HC_\idot(A) \cong
HH_\idot(A^\tw)[v]$.
\end{lemma}

\proof{} By Lemma~\ref{dag.dege}, the conjugate spectral sequence
for $W_1HC_\idot(A)$ degenerates up to term $E^{p-2}$. By the
assumption \ref{p.large} of Theorem~\ref{char.p.main}, it
degenerates in all the following terms for dimension reasons.
\endproof

We now recall that by definition, we have an exact triangle
\begin{equation}\label{defekt}
\begin{CD}
W_1HC_\idot(A) @>>> HC_\idot(A) @>>> H^\hdot(\langle Z,J
\rangle,HC_\idot(\overline{A}_{\#})) @>>>
\end{CD}
\end{equation}
of complexes of $k$-vector spaces, where the cyclic object
$\overline{A}_{\#} \in \Fun(\Lambda,\C)$ is the quotient
$$
\overline{A}_{\#} = \pi_!i_p^*A_{\#}/A^\tw_{\#}.
$$

\begin{lemma}\label{bounded.dg}
In the assumptions of Theorem~\ref{char.p.main}, the Hochschild
homology sheaf $HH_l(\overline{A}_{\#}) \in \C$ is trivial whenever
$l > mp$.
\end{lemma}

\proof{} As in the proof of Lemma~\ref{otimesp}, we may assume that
$\langle Z,J \rangle$ is trivial, so that $\C$ is the category of
$k$-vector spaces. Moreover, the homology of $A^\tw_{\#}$ is bounded
by Theorem~\ref{char.p.main}~\ref{reso}, so that we might just as
well prove our statement for $\pi_!i_p^*A_{\#}$ instead of its
quotient $\overline{A}_{\#}$.

For any $A$-bimodule $M$, denote by $M_{\bekar}$
the cokernel of the commutator map $A \otimes M \to M$ --
equivalently, we have $M_{\bekar} = M \otimes_{A \otimes A^{opp}}
A$. Moreover, denote 
$$
M_{p\bekar} = \left(M \otimes_A M \otimes_A \dots \otimes_A
M\right)_{\bekar},
$$
where we have $p$ multiples on the right-hand side. The reader will
check easily that this is construction is actually cyclically
symmetric in all $p$ multiples (but since paper is two-dimensional,
we cannot represent this symmetry in convenient notation). Moreover,
$M \mapsto M_{p\bekar}$ is obviously functorial in $M$.

Recall that by the Dold-Thom Theorem \cite{DT}, for any abelian
category $\C$ the category $\Fun(\Delta^o,\C)$ of simplicial objects
in $\C$ is equivalent to the category of complexes in $\C$
concentrated in non-positive degree. For any $M \in
\Fun(\Delta^o,\C)$, the corresponding complex represents the
homology $H_\idot(\Delta^o,M)$. If the category $\C$ is equipped
with a symmetric tensor product, then the category
$\Fun(\Delta^o,\C)$ also has a natural symmetric tensor structure,
and this structure is homotopy-invariant in the following strong
sense. For any group $G$, denote by $\C[G]$ the category of objects
in $\C$ equipped with an action of $G$.  For any $M \in
\Fun(\Delta^o,\C)$, the $p$-fold tensor product $M^{\otimes p}$
equipped with the transposition action of the symmetric group $S_n$
is naturally an object in $\Fun(\Delta^o,\C[S_n])$ -- or, by
Dold-Thom equivalence, a complex in $\C[S_n]$. Then if $M_0,M_1 \in
\Fun(\Delta^o,\C)$ are homotopy-equivalent -- that is, the
corresponding complexes are equivalent up to chain homotopy of
complexes -- then the complexes $M_o^{\otimes p}$, $M_1^{\otimes p}$
are homotopy-equivalent in $\C[S_n]$.

If the category $\C$ is tensor but not symmetric, there is no
natural group action on tensor powers $M^{\otimes p}$. However, if
$\C=A\bimod$ is the category of $A$-bimodules, then for any
simplicial $A$-bimodule $M \in \Fun(\Delta^o,\C)$ we can form
$M_{p\bekar} \in \Fun(\Delta^o,\C[\Z/p\Z])$, and this is also
homotopy-invariant in the above sense.

Take now a flat resolution $P_\idot$ of the diagonal $A$-bimodule
$A$ of length $m$ whose existence is assumed by
Theorem~\ref{char.p.main}~\ref{reso}, and treat it as a simplicial
$A$-bimodule $P_\idot \in \Fun(\Delta^o,A\bimod)$. Moreover, let
$P'_\idot \in \Fun(\Delta^o,A\bimod)$ be the standard simplicial
bar-resolution of $A$. Then by definition of the cyclic object
$A_{\#} \in \Fun(\Lambda,k\mmod)$, we have
$$
j^*i_p^*A_{\#} \cong \left(P'_\idot\right)_{p\bekar} \in
\Fun(\Delta^o,k[\Z/p\Z]\mmod) = \Fun(\Delta^o \times \ppt_p,k\mmod),
$$
and since $P_\idot$ must be homotopy-equivalent to $P'_\idot$, this
is homotopy-equivalent to $(P_\idot)_{p\bekar}$ as a complex of
$k[\Z/p\Z]$-modules. We conclude that
$$
j^*\pi_!i_p^*A_{\#} \cong
H_0\left(\Z/p\Z,\left(P'_\idot\right)_{p\bekar}\right) 
$$
is homotopy equivalent to $H_0(\Z/p\Z,(P_\idot)_{p\bekar}) \in
\Fun(\Delta^o,k\mmod)$. Applying the Dold-Thom equivalence, we see
that the latter complex is manifestly trivial in degrees $> mp$.
\endproof

\proof[Proof of Theorem~\ref{char.p.main}.] Consider the exact
triangle \eqref{defekt}. We have
$$
HC_\idot(A) \cong H^\hdot(\langle Z,J
\rangle,H_\idot(\Lambda_p,i_p^*A_{\#})) \cong H^\hdot(\langle Z,J
\rangle,H_\idot(\Lambda,\Ll^\hdot\pi_!i_p^*A_{\#})),
$$
and by Lemma~\ref{HS}, the generator $u \in H^2(\Lambda,k)$ of the
cohomology algebra $H^\hdot(\Lambda,k)$ acts trivially on the
right-hand side. Therefore the cokernel of the natural map
$W_1HC_\idot(A) \to HC_\idot(A)$ lies inside the subspace in
$H^\hdot(\langle Z,J \rangle,HC_\idot(\overline{A}_{\#}))$
annihilated by $u$. This subspace is a quotient of $H^\hdot(\langle
Z,J \rangle,HH_\idot(\overline{A}_{\#}))$, and by
Lemma~\ref{bounded.dg}, the latter space is trivial in degrees $\lle
0$. We conclude that the map $W_1HC_\idot(A) \to HC_\idot(A)$ is
trivial in degrees $\lle 0$. By Lemma~\ref{conj.splits}, this implies
that
$$
\dim HC_l(A) \leq \bigoplus_m \dim HH_{l+2m}(A^\tw) = \bigoplus_m
\dim HH_{l+2m}(A)
$$
for $l \lle 0$ (the sum on the right-hand side is over all integers
$m$, and it is bounded and finite by our assumptions). By the
standard criterion \cite{de}, this implies that the
Hochschild-to-cyclic spectral sequence $HH_\idot(A)[v] \Rightarrow
HC_\idot(A)$ degenerates in degrees $\lle 0$. Since the sequence is
periodic, it degenerates everywhere.
\endproof

This finishes the positive characteristic case. The characteristic
$0$ statement would follow from this by a standard and well-known
procedure, but there is a problem: interesting non-commutative
algebras are often not Noetherian. Therefore it is unclear how to
``spread out'' an algebra given over a field of characteristic $0$
so that it acquires fibers over fields of positive
characteristic. In this paper, we have decided to settle for the
following comporimise: we give a statement ``in mixed
characteristic'' -- that is, we assume that the algebra is already
defined and flat over a $\Z$-algebra of finite type -- and we give a
second statement which shows that such a such a $\Z$-model exists in
the Noetherian situation. This is enough for some applications but
certainly not for all of them. This is a technical limitation, since
the finiteness assumptions we impose in order to have the
degeneration should also be sufficient to construct a $\Z$-model. In
a separate paper, we will handle this problem by using the
$A_\infty$-methods.

\begin{theorem}\label{model.main}
Assume given an associative unital algebra $A \in \Shv(Z,J)$ flat
over a intergal domain $O$ of finite type over $\Z$. Moreover,
assume that
\begin{enumerate}
\item the diagonal $A$-bimodule $A$ admits a finite flat resolution,
\item for some integer $m > 0$ and every sheaf $E \in \Shv(Z,J)$ of
  $K$-vector spaces on $\langle Z,J \rangle$, the cohomology groups
  $H^l(\langle Z,J \rangle,E)$ are trivial whenever $l > m$, and
\item for any integer $l$, the Hochschild homology group $HH_l(A)$
  is a finitely generated $O$-module.
\end{enumerate}
Then the Hochschild-to-cyclic spectral sequence $HH_\idot(A)[v]
\Rightarrow HC_\idot(A)$ degenerates.
\end{theorem}

\proof{} Denote by $n$ the length of the resolution in
\thetag{i}. It suffices to show that the differentials in the
spectral sequence vanish after reduction modulo all maximal ideals
$\m \subset O$, and moreover, it suffices to consider all the
maximal ideals outside of a closed subset in $\Spec O$. Therefore we
may assume that $\cchar O/\m > n + m$. Since $O$ is of finite type
over $\Z$, the field $O/\m$ is perfect, and we are done by
Theorem~\ref{char.p.main}.
\endproof

\begin{theorem}\label{char.0.main}
Assume given a field $K$ of characteristic $0$, and assume given an
associative unital Noetherian $K$-algebra $A \in \Shv(Z,J)$ which is
generated by a finite number of local sections. Moreover, assume
that
\begin{enumerate}
\item the diagonal $A$-bimodule $A$ admits a finite flat resolution,
\item for some integer $m > 0$ and every sheaf $E \in \Shv(Z,J)$ of
  $K$-vector spaces on $\langle Z,J \rangle$, the cohomology groups
  $H^l(\langle Z,J \rangle,E)$ are trivial whenever $l > m$, and
\item for any integer $l$, the Hochschild homology group $HH_l(A)$
  is a finite-dimensional $K$-vecotr space.
\end{enumerate}
Then the Hochschild-to-cyclic spectral sequence $HH_\idot(A)[v]
\Rightarrow HC_\idot(A)$ degenerates.
\end{theorem}

\proof{} Since $A$ is generated by a finite number of local
sections, there exists a subalgebra $O \subset K$ of finite type
over $\Z$ such that there exists an $O$-algebra $A_O \in \Shv(Z,J)$
with $A \cong A_O \otimes_O K$. Moreover, we may assume that the
finite flat resolution of the diagonal bimodule $A$ is defined over
$O$ (as a complex, possibly not acyclic). Since $A$ is Noetherian,
the homology of this complex consists of finitely generated torsion
$O$-modules; localizing $O$, we may assume that the complex is also
a resolution over $O$. Finally, again since $O$ is Noetherian, we
may assume after localizing $O$ that the resolution is flat, and
that $A_O$ itself is flat over $O$. Now $A_O$ satisfies all the
assumptions of Theorem~\ref{model.main}.  \endproof

\bigskip

\noindent
{\sc
Steklov Math Institute \& ITEP\\
Moscow, USSR}

\bigskip

\noindent
{\em E-mail address\/}: {\tt kaledin@mccme.ru}

\end{document}